\documentclass{article}

\usepackage{amssymb,epsf}
\usepackage{amsfonts,amssymb}
\usepackage{amsmath,bbm}
\usepackage{theorem}

\newtheorem{tht}{Theorem}%[chapter]
\newtheorem{thd}{Definition}%[chapter]
\newtheorem{thl}[tht]{Lemma}
\newtheorem{thp}[tht]{Proposition}
\newtheorem{thc}[tht]{Corollary}

\theoremstyle{plain} {\theorembodyfont{\rmfamily}
\newtheorem{thex}{Example}%[chapter]
\newcommand{\dK}{\mathbb{K}}
\newcommand{\mn}{\medskip\noindent}

\newcommand{\cP}{{\mathcal{P}}}

\newcommand{\cW}{{\mathcal{W}}}
\newcommand{\cS}{{\mathcal{S}}}

\newcommand{\cN}{{\mathcal{N}}}
\newcommand{\cC}{{\mathcal{C}}}
\newcommand{\cD}{{\mathcal{D}}}
\newcommand{\cJ}{{\mathcal{J}}}

\newcommand{\Hh}{{\mathcal{H}}}
\newcommand{\cG}{{\mathcal{G}}}
\newcommand{\cK}{{\mathcal{K}}}

\newcommand{\cB}{{\mathcal{B}}}

\newcommand{\cO}{{\mathcal{O}}}
\newcommand{\cZ}{{\mathcal{Z}}}
\newcommand{\cR}{{\mathcal{R}}}

\newcommand{\cI}{{\mathcal{I}}}
\newcommand{\cX}{{\mathcal{X}}}
\newcommand{\cT}{{\mathcal{T}}}
\newcommand{\cA}{{\mathcal{A}}}
\newcommand{\cE}{{\mathcal{E}}}
\newcommand{\cM}{{\mathcal{M}}}
\newcommand{\nein}[1]{}

\newcommand{\dR}{\mathbb{R}}
\newcommand{\dN}{\mathbb{N}}
\newcommand{\dC}{\mathbb{C}}

\newcommand{\ov}{\overline}

\newcommand{\tint}{{\textstyle\int}}

\newcommand{\ii}{\mathrm{i}}

\newcommand{\aabs}[1]{\left\| #1\right\|}

\hyphenation{ho-mo-mor-phism in-var-i-ant rep-re-sen-ta-tion
            al-ge-bra sub-al-ge-bra}

\newcommand{\uhr}{{\upharpoonright}}
\renewcommand{\lceil}{\uhr}

\begin{document}

\title{NONCOMMUTATIVE REAL ALGEBRAIC GEOMETRY -\\ SOME BASIC CONCEPTS AND FIRST IDEAS}
\author{KONRAD SCHM\"UDGEN\thanks{Fakult\"at f\"ur Mathematik und Informatik, Universit\"at Leipzig, Germany}}

\maketitle

\begin{abstract}
We propose and discuss how basic notions (quadratic modules,
positive elements, semialgebraic sets, Archimedean orderings) and
results (Positivstellens\"atze) from real algebraic geometry can be
generalized to noncommutative $\ast$-algebras. A version of
Stengle's Positivstellensatz for $n \times n$ matrices of real
polynomials is proved.
\end{abstract}

%\keywords{Noncommutative real algebraic geometry, quadratic module,
%sum of squares, $\ast$-representation, positivity, positive
%semidefinite matrices}

%\subjclass{Primary: 13J30, 47L60, 15A48 Secondary: 11E25}

\section{Introduction}
In recent years various versions of noncommutative
Positivstellens\"atze have been proved, about free polynomial
algebras by J.W. Helton and his coworkers \cite{hel},\cite{hm},
\cite{hmp} and about the Weyl algebra \cite{sweyl} and enveloping
algebras of Lie algebras \cite{senv} by the author. These results
can be considered as very first steps towards a new mathematical
field that might be called {\it noncommutative real algebraic
geometry}.

The aim of this paper is
%to propose and
to discuss how some basic
concepts and results from real algebraic geometry should look like
in the  noncommutative setting. This article unifies a series of
talks I gave during the past five years at various conferences
(Pisa, Marseille, Palo Alto, Banff) and at other places. It should
be emphasized that it represents the authors personal view and
ideas on this topic. These concepts and ideas exist and will be
presented at different levels of exactness and acceptance by the
community. Some of them (quadratic modules, positivity by
representations, Archimedean orderings) are more or less clear and
accepted. The definition of a semialgebraic set seems to be
natural as well. Possible formulations of Artin's theorem in the
noncommutative case are at a preliminary stage and will be become
clearer when more results are known, while others such as
the definition of the preorder or a noncommutative formulation of
Stengle's theorem require more
research before satisfying formulations can be given.

In section 2 we collect some general definitions and notations
which are used throughout this paper. In section 3 we set up basic
axioms, concepts and examples for noncommutative real algebraic
geometry. Roughly speaking, by passing to the noncommutative case
the polynomial algebra and the points of $\dR^d$ are replaced by a
finitely generated $\ast$-algebra and a distinguished family of
irreducible $\ast$-representations. In section 4 we investigate
and discuss possible formulations of Artin's theorem and Stengle's
theorem in the noncommutative case. As the main new result of this
paper we obtain a version of Stengle's theorem for the algebra of
$n \times n$ matrices with polynomial entries. Section 5 is
devoted to $\ast$-algebras with Archimedean quadratic modules. We
derive some properties and abstract Positivstellens\"atze for such
$\ast$-algebras and develop a variety of examples. In
section 6 we show how Pre-Hilbert $\ast$- bimodules can be used to
transport quadratic modules from one algebra to another. This
applies nicely to algebras of matrices and it might
have further applications.

I thank Y. Savchuk for helpful discussions on the
subject of this paper.

\section{Definitions and Notations}

{\it Throughout this article $\cA$ denotes a real or complex
unital $\ast$-algebra.} By a $\ast$-{\it algebra} we mean an
algebra $\cA$ over the field $\dK=\dR$ or $\dK=\dC$ equipped with
a mapping $a\to a^\ast$ of $\cA$ into itself, called the {\it
involution} of $\cA$, such that $(\lambda a+\mu b)^\ast =
\bar{\lambda}a^\ast+ \bar{\mu} b^\ast, (ab)^\ast = b^\ast a^\ast$
and $(a^\ast)^\ast=a$ for $a,b\in \cA$ and $\lambda, \mu\in \dK$.
The unit element of $\cA$ is denoted by $1$ and $\cA_h :=\{a\in\cA:
a=a^\ast\}$ is the set of {\it hermitian elements}
of $\cA$.

As usual, $\dR[t]{=}\dR[t_1,{\dots},t_d]$ resp.
$\dC[t]{=}\dC[t_1,{\dots},t_d]$ are the $\ast$-algebras of real
resp. complex polynomials in $d$ commuting hermitian
indeterminates $t_1,{\dots},t_d$. Set $\cN(p)=\{s \in
\dC^d:p(s)=0\}$ for $p \in \dC[t]$. Let $\cM_{k,n}(R)$ denote the ${k\times
n}$-matrices over a ring $R$ and set $\cM_n(R):=\cM_{n,n}(R)$.

If $a$ is an operator on a Hilbert space, we denote by
$\cD(a)$ its {\it domain}, by $\cR(a)$ its {\it range}, by
$\cN(a)$ its kernel, by $\bar{a}$ its {\it closure} and by $a^\ast$
its {\it adjoint} (if they exist). A subset $\cE$
of $\cD(a)$ is called a {\it core} for $a$ if for
each $\varphi \in \cD(a)$ there is a sequence of vectors
$\varphi_n \in \cE$ such that $\varphi_n \to \varphi$ and
$a\varphi_n \to a\varphi$.

We now turn to some notions on $\ast$-representations, see
\cite{unb} for a treatment of this subject. Let $\cD$ be a
pre-Hilbert space with scalar product $\langle\cdot,\cdot\rangle$.
A $\ast$-{\it representation} of $\cA$  on
$\cD(\pi){:=}\cD$ is an algebra homomorphism $\pi$ of $\cA$ into
the algebra of linear operators mapping $\cD$ into itself
such that $\pi (1)\varphi=\varphi$ and $\langle
\pi(a)\varphi,\psi\rangle =\langle \varphi,\pi(a^\ast)\psi\rangle$
for all $\varphi,\psi\in\cD$ and $a\in \cA$. Two
$\ast$-representation $\pi_1$ and $\pi_2$ are {\it
(unitarily) equivalent} if there exists an isometric linear
mapping $U$ of $\cD(\pi_1)$ onto $\cD(\pi_2)$ such that $\pi_2 (a)
= U\pi_1(a)U^{-1}$ for $a\in \cA$. A $\ast$-representation
$\pi$ is called {\it irreducible} if any decomposition of
$\cD(\pi)$ as an orthogonal sum of subspaces $\cD_1$ and $\cD_2$
such that $\pi (a) \cD_1 \subseteq \cD_1$ and $\pi (a) \cD_2
\subseteq\cD_2$ for all $a\in \cA$ implies that $\cD_1=\{0\}$ or
$\cD_2=\{0\}$.

A {\it state} of $\cA$ is a linear functional $f$ on $\cA$ such
that $f(1)=1$ and $f(a^\ast a) \geq 0$ for all $a \in \cA$. A
state $f$ of $\cA$ is called {\it pure} if each state $g$
satisfying $g(a^\ast a)\leq f(a^\ast a)$ for all $a \in \cA$ is a
multiple of $f$. If $\pi$ is a $\ast$-representation of $\cA$ and
$\varphi\in \cD(\pi)$ is a unit vector, then $f(\cdot):=\langle
\pi(\cdot)\varphi, \varphi \rangle $ is a state of $\cA$. These
states are called {\it vector states} of $\pi$. Each state arises
in this manner. That is, for each state $f$ on $\cA$ there exists
a distinguished $\ast$-representation of $\cA$, called the {\it
GNS-representation} of $f$ and denoted by $\pi_f$ (see \cite{unb},
8.6), and a vector $\varphi \in \cD(\pi_f)$ such that $\cD(\pi_f)=
\pi_f(\cA)\varphi_f$ and
\begin{align}\label{gns}
f(a) = \langle \pi_f(a) \varphi_f,\varphi_f \rangle ~{\rm for}~ a \in \cA.
\end{align}
We recall two standard definitions from real algebraic
geometry (see e.g. \cite{bcr}, \cite{pd} or \cite{mm}). Suppose
$\cB$ is a commutative unital algebra over $\dK=\dR,\dC$. Let $\hat{\cB}$ denote
the set of all homomorphisms of $\cB$ into $\dK$. We write
$f(s):=s(f)$ for $f\in \cB$ and $s \in \hat{\cB}$. If
$\cB=\dR[t_1,\dots,t_d]$, then each element of $\hat{\cB}$ is given by the
evaluation at some point of $\dR^d$, that is, $\hat{\cB}\cong
\dR^d$.

If $f{=}(f_1,\cdots,f_k)$ is a k-tuple of elements $f_j \in \cB$,
the {\it basic closed semialgebraic set} $\cK_f$ and the {\it
preorder} $\cT_f$ associated with $f$ are defined by
\begin{align}\label{kf}
\cK_f &=\{ s \in \hat{\cB}: f_1(s) \geq 0,\cdots, f_r(s)\geq 0
\},\\\label{tf}
\cT_f &=\{~\textstyle{\sum_{\varepsilon_i \in
\{0,1\}}}~\textstyle{ \sum_{l=1}^r}~ f_1^{\varepsilon_1}\cdots
f_k^{\varepsilon_k} g_l^2 ~;~ g_l \in \cB,~r \in \dN \}.
\end{align}

\section{Basic Concepts of Noncommutative Real Algebraic Geometry}\label{s1}
\subsection{Two Main Ingredients of Noncommutative Real Algebraic Geometry}\label{ss1}

The first main ingredient of real algebraic geometry is the
algebra $\dR[t_1,{\dots},t_d]$ of polynomials or the coordinate
algebra $\dR[V]$ of a real algebraic variety. Its counter-part in
the noncommutative case is a
$$\centerline{$\bullet$ {\it finitely generated real or complex unital
$\ast$-algebra $\cA$}.}
$$
In real algebraic geometry elements of the algebras
$\dR[t_1,{\dots}, t_d]$ or $\dR[V]$ are evaluated at points of
$\dR^d$ or of $V$. As the noncommutative substitute
for the set of point evaluations we assume that we have a given distinguished
$$
\centerline{$\bullet$ {\it family $\cR$ of equivalence classes of
irreducible $\ast$-representations of $\cA$.}}
$$
Elements of $\cR$ can be interpreted as "points" of a "noncommutative space".

One could also take a {\it family of pure states on $\cA$} instead
of representations. The GNS-representation of a pure state
is irreducible (\cite{unb},8.6.8). The converse is valid for bounded
representations on a Hilbert space.

If $\pi$ is a $\ast$-representation of $\cA$ and $a\in
\cA_h$, we write
%$$
%\pi (a)\ge 0~{\rm if~and~only~if}~ \langle \pi (a) %\varphi,\varphi\rangle \ge 0~{\rm of~all}~ \varphi \in\cD(\pi).
%$$
$$
\pi (a)\ge 0~~if~and~only~if~~ \langle \pi (a) \varphi,\varphi\rangle \ge 0~~ for~all~~ \varphi \in\cD(\pi).
$$
This may be considered as a generalization of the positivity
$f(t)\ge 0$ of the point evaluation of $f\in\dR[t]$ at $t\in\dR^d$.

Let us collect a number of important examples.

\begin{thex}\label{pol1}{\it Commutative polynomial algebras}\\
The $\ast$-algebra for "ordinary" real algebraic geometry is the
real $\ast$-algebra $\cA=\dR [t_1,{\dots},t_d]$ with trivial
involution $a^\ast{=}a$. One can also take the complex
$\ast$-al\-geb\-ra $\cA=\dC [t_1,{\dots},t_d]$ with involution
defined by $p^\ast(t){:=}\sum_\alpha \overline{c}_\alpha~ t^\alpha$
for $p(t){=}\sum_\alpha c_\alpha t^\alpha $.

Let $\cR \cong \dR^d$ be the set of point evaluation, that is,
$\cR=\{\pi_t : t\in \dR^d\}$, where $\pi_t (p)=p(t)$ for $ p\in \cA$.
 \hfill $\circ$
\end{thex}

\begin{thex}\label{weyl1} {\it Weyl algebras}\\
Let $d \in \dN$. The Weyl algebra $\cW(d)$ is the unital
$\ast$-algebra with generators $a_1,\dots,a_d$,
$a_{-1},\dots,a_{-d}$,~ defining relations
$$
a_ka_{-k}-a_{-k}a_k =1 ~~{\rm and}~~ a_ka_l =a_la_k ~~ {\rm if}~~ k\neq l,
$$
and involution defined by~ $(a_k)^\ast {=}a_{-k}$ for
$k{=}1,\dots,d$.

For the Weyl algebra $\cW(d)$ the set $\cR$ consists of a single
element $\pi_0$, the {\it Bargmann-Fock representation}.
It is described by the actions of generators on an orthonormal basis
$\{e_n;n \in \dN_0^d \}$ of the Hilbert space given by
$$
\pi_0(a_k)e_n = n_k^{1/2} e_{n-1_k},~ \pi_0(a_{-k})e_n =(n_k+1)^{1/2} e_{n+1_k}
$$
for $k{=}1,{\dots},d$ and $n{=}(n_1,\dots,n_d) {\in} \dN_0^d$. Here
$1_k$ is the d-tuple with $1$ at the k-th place and $0$
otherwise and $e_{n-1_k}{:=}0$ when $n_k=0$. The domain $\cD(\pi_0)$
consist of all sums $\sum \varphi_n e_n$ such that $\sum
n_1^r\dots n_d^r |\varphi_n|^2 < \infty$ for all $r \in \dN$. \hfill
$\circ$
\end{thex}

\begin{thex}\label{env1} {\it Enveloping algebras}\\
Let $\cE(\cG)$ be the complex universal enveloping algebra of a
finite dimensional real Lie algebra $\cG$. Then $\cE(\cG)$ is a
$\ast$--algebra with involution given by $x^\ast {=}-x$ for $x \in
\cG$.

There is a simply connected Lie group $G$ having $\cG$ as its Lie
algebra. Let $\hat{G}$ denote the unitary dual of $G$, that is,
$\hat{G}$ is the set of equivalence classes of irreducible unitary
representations of $G$. For each $\alpha \in \hat{G}$ we fix a
representation $U_\alpha$ of the class $\alpha$. Each
representation $U_\alpha$ of $G$ gives rise to an irreducible
$\ast$-representation $dU_\alpha$ of the $\ast$-algebra
$\cE(\cG)$, see e.g. \cite{unb}, 10.1, for details. The
irreducibility of $dU_\alpha$ follows from \cite{unb}, 10.2.18.

As $\cR$ we take the family $\{dU_\alpha; \alpha {\in} \hat{G} \}$ of
$\ast$-representations of $\cE(\cG)$. \hfill $\circ$
\end{thex}

\begin{thex}\label{free1} {\it Free polynomial algebras}\\
For $d \in \dN$, let $\cA=\dC\left< t_1,\dots,t_d \right>$ be the
free unital complex algebra with $d$ generators $t_1,\dots,t_d$. It
is a $\ast$-algebra with involution determined by $t_j^\ast =t_j$,
$j{=}1,\dots,d$.

Let $\cR$ be the equivalence classes of all irreducible
representations by bounded operators on a Hilbert space. We may also
take the families $\cR_1$ of equivalence classes of $\pi$ in $\cR$
which act on a {\it fixed} (sufficiently large) Hilbert space or
$\cR_2$ of equivalence classes of finite dimensional representations
in $\cR$. Since the $\ast$-al\-geb\-ra is free, {\it each} $d$-tuple
of bounded selfadjoint operators $T_j$ on a Hilbert space $\Hh$
defines a $\ast$-representation $\pi$ on $\Hh$ by $\pi(t_j)=T_j$.

Often it is convenient to use the $\ast$-algebra
$\cA_0{=}\dC\left<z_1,{\dots},z_d,w_1,{\dots},w_d\right>$ with
involution given by $z_j^\ast:=w_j$. Clearly, $\cA_0$ is
$\ast$-isomorphic to
%$\dC\left<t_1,\dots,t_{2d}\right>$
$\cA$ with a
$\ast$-isomorphism determined by $z_j \to t_j+\ii t_{d+j}$,
$j{=}1,{\dots,}d$. \hfill $\circ$
\end{thex}

\begin{thex}\label{matcomm} {\it Matrix algebras over commutative $\ast$-algebras}\\
Suppose that $\cB$ is a commutative unital $\ast$-algebra. Let $\cA$
be the matrix $\ast$-algebra $\cM_n(\cB)$ with involution
$(b_{ij})^\ast = (b_{ji}^\ast)$ and $\cR$ the set $\{\rho_s; s \in
\dR^d\}$ of irreducible $\ast$-representations $\rho_s{:}A \to
A(s)$, where the matrix $A(s)$ acts as linear operator on the
Hilbert space $\dK^d$ in the usual way.

{\it $\ast$-Subalgebras of matrix $\ast$-algebras} $\cM_n(\cB)$
provide a large class of interesting $\ast$-algebras for
noncommutative real algebraic geometry. Example \ref{yurii} below is
one of such examples.  More can be found in the book \cite{os}.
\end{thex}
\mn How a possible new theory as noncommutative real algebraic
geometry will evolve in future depends essentially on what will be
considered as typical examples and fundamental problems.
Positivstellens\"atze should be one of the basic problems to be
studied. $\ast$-Subalgebras of matrix algebras over commutative
$\ast$-algebras (see Example \ref{matcomm}) will lead to a theory
that is closest to real algebraic geometry. They might be studied
first. Weyl algebras and enveloping algebras (and some algebras from
subsection \ref{exarch}) are interesting but  challenging classes of
examples. It is likely that a theory based on these examples will be
very different from the classical theory.

\subsection{Quadratic Modules and Orderings}
\begin{thd}
A {\rm quadratic module} of $\cA$ is a subset $\cC$ of  $\cA_h$ such that
\begin{eqnarray}\label{QM1}
&1 \in \cC,~~\cC + \cC \subseteq \cC,~~ \dR_+{ \cdot}\cC \subseteq \cC,\\
\label{QM2}
&b^\ast \cC b \in \cC~{\rm for~all}~b\in \cA.
\end{eqnarray}
\end{thd}
Quadratic modules are  important in theory of $\ast$-algebras where
they have been called {\it $m$-admissible wedges} (\cite{unb}, p.
22). Following the  terminology from real algebraic geometry we
prefer to use the name "quadratic module".

Each quadratic module gives an ordering $\preceq$ on the real vector
space $\cA_h$ by defining $a \preceq b$ (and likewise $b \succeq a$)
if and only if $a-b \in \cC$.

All elements $a^\ast a$, where $a \in \cA$, are called  {\it squares} of $\cA$. The wedge
$$
\sum \cA^2 := \left\{ \sum^n_{j=1} a_j^\ast a_j;~~ a_1,{\dots},a_n\in \cA, n\in \dN \right\}
$$
of finite sums of squares is obviously the smallest quadratic module of $\cA$.

If $\cS$ is a family of $\ast$-representations of $\cA$, then
$$
\cA(\cS)_+ := \{ a\in \cA_h: \pi(a) \geq 0~~{\rm for~all} ~~\pi \in \cS\}
$$
is a quadratic module of $\cA$. The interplay between quadratic
modules which are defined in algebraic terms (such as $\sum
\cA^2$)
and those which are defined by means of $\ast$-representations (such as $\cA(\cS)_+$)
is one of the most interesting challenge for the theory.

The following polarization identities are useful. For
 $a,x,y {\in} \cA$, we have
\begin{align}\label{id1}
4x^\ast ay = ~~&(x+y)^\ast a (x+y)-(x-y)^\ast a (x-y)\\\nonumber&-\ii(x+ \ii y)^\ast
a (x + \ii y)+\ii(x-\ii y)^\ast a (x-\ii y),\\
2 (x^\ast ay + y^\ast ax)~=~&(x+y)^\ast a(x+y)-(x-y)^\ast a (x-y).\label{id2}
\end{align}
From (\ref{id2}), applied with $a{=}y{=}1$, and (\ref{id1})  we easily conclude that
\begin{align}\label{id3}
A_h = \cC - \cC,~~ \cA=(\cC-\cC)+\ii (\cC - \cC).
\end{align}
for any quadratic module $\cC$. Of course, for (\ref{id1}) and for
the second equalitiy of (\ref{id3}) one has to assume that $\cA$ is
a {\it complex} $\ast$-algebra.

A quadratic module $\cC$ is called {\it proper} if $\cC \ne \cA_h$.
By (\ref{id2}),  $\cC$ is proper iff $-1$ is not in $\cC$. A proper
quadratic module $\cC$ of $\cA$ is called {\it maximal} if there is
no proper quadratic module $\tilde{\cC}$ of $\cA$ such that $\cC
\subseteq \tilde{\cC}$ and $\cC\ne \tilde{\cC}$.

If $\cC$ is a maximal proper quadratic module of a commutative
unital ring $A$, then $\cC\cap (-\cC)$ is a prime ideal  and $\cC
\cup (-\cC)=A$. In the noncommutative case the second assertion is
not true, for the first we have the following theorem due to J.
Cimpric \cite{cim1}.

\begin{tht}\label{cimpric}
Suppose  $\cC$ is a quadratic module of a complex $\ast$-algebra
$\cA$. Let $\cC^0 := \cC \cap(-\cC)$ and $\cI_\cC:= \cC^0 + \ii
\cC^0.$
{\rm (i)} $\cI_\cC$ is a two-sided $\ast$-ideal of $\cA$.\\
{\rm (ii)} If $\cC$ is a maximal proper quadratic module, $\cI_\cC$ is a prime ideal and
$$
\cI_\cC =\{ a\in\cA: a xx^\ast a^\ast \in\cC^0~{\rm for~all}~x\in\cA\}.$$
\end{tht}
{\bf Proof.} (i) Clearly, $\cI_\cC$ is $\ast$-invariant
and  $\cC^0$ is a real vector subspace. If $a\in \cC^0$ and
$x\in \cA$, then $(x^\ast {+} \ii^ky)^\ast a (y^\ast {+} \ii^k y)\in \cC^0$
for $k{=}0,1,2,3$ by (\ref{QM2}) and hence $4xay\in\cC^0 + \ii  \cC^0 {=}
\cI_\cC$ by (\ref{id1}). Thus
$\cA\cdot\cC^0\cdot\cA\subseteq\cJ_\cC$ and hence
$\cA\cdot \cI_{\cC}\cdot \cA \subseteq \cI_\cC.$\\
(ii) \cite{cim1}, Theorem 1 and Remark on p.5.\hfill $\Box$

\subsection{Noncommutative Semialgebraic Sets}\label{ncsas}
Let $\cR$ be a family of (equivalence classes) of $\ast$-representations of $\cA$.
\begin{thd}\label{def1}
A subset $\cK$ of $\cR$ is {\rm semialgebraic} if it is a finite Boolean
combination (that is, using unions, intersections and complements) of sets $\{ \pi
\in \cR : \pi (f) \ge 0 \}$ for $f \in \cA_h$. It is {\rm
algebraic} if there is a finite subset $f {=} \{ f_1, {\dots}, f_k \}$ of $\cA$ such that
$\cK =\cZ(f){:= }\{ \pi \in \cR : \pi (f_1) {=}0, {\cdots}, \pi (f_k){=} 0 \}$.
\end{thd}
Let $f=(f_1,{\dots},f_k)$ be a k-tuple of elements of $\cA_h$, where
$f_1=1$. We define the {\it basic closed semialgebraic set}
associated with $f$ by
\begin{align}\label{bas}
\cK(f)= \{\pi\in\cR:\pi(f_1)\ge 0,{\dots}, \pi (f_k) \ge 0\}
\end{align}
and the associated wedges
\begin{eqnarray}\label{pf1}
&\cP(f)= \{ a\in \cA_h:\pi (a)\ge  0~{\rm
for~all}~\pi\in\cK(f)\},\\\label{tf1} &\cC(f)=\{~ \textstyle{
\sum^s_{j=1}} \textstyle{\sum^k_{l=1}}~ a^\ast_{jl} f_l
a_{jl}:a_{jl}\in\cA,~s \in \dN\}.
\end{eqnarray}
Then $\cP(f)$ and $\cC(f)$ are quadratic modules such that $\cC(f)
\subseteq \cP(f)$ and $\cC(f)$ is the smallest quadratic
module that contains all elements $f_1,\dots,f_k$. In general we
cannot add mixed products $f_jf_l$ to the wedge $\cC(f)$, because
$f_jf_l$ is not hermitian if $f_j$ and $f_l$ do not commute.

\begin{thex}\label{pol2}{\it Commutative polynomial algebras}\\
If $\cA$ and $\cR$ are as in Example \ref{pol1}, semialgebraic sets,
algebraic sets, and basic closed semialgebraic sets according to the
preceding definitions are just the ordinary ones in real algebraic
geometry \cite{bcr},\cite{pd}\cite{mm} and $\cP(f)$ is the wedge of
nonnegative real polynomials on $\cK(f)$. Since $\cC(f)$ is in
general not closed under multiplication, it is not a preorder. If we
replace $f$ by the tuple $\tilde{f}$ of all products $f_{i_1}\dots
f_{i_r}$, where $1\leq i_1 <i_2<\dots < i_r \leq k$, then
$\cK(f)=\cK({\tilde{f}})$ and $\cC({\tilde{f}})$ is the usual
preorder $\cT_f$. \hfill $\circ$
\end{thex}

\begin{thex}\label{free2} {\it Free polynomial algebras}\\
Let $\cA=\dC\left<t_1,{\dots},t_d\right>$ and
$\cA_0=\dC\left<z_1,{\dots},z_d,w_1,{\dots,},w_d \right>$ be the
$\ast$-algebras from Example \ref{free1} and let $\cR$ be the family
of all bounded $\ast$-representations of $\cA$ resp. $\cA_0$ on a
separable Hilbert space. Then the basic semialgebraic set $\cK(f)$
for $\cA$ defined by (\ref{bas}) corresponds precisely to the {\it
positivity domain} $\cD_f$ of $f$ according to J. Helton and S.
McCullough \cite{hm}. For the polynomial $f(z_1,\dots,z_d):=
z_1^\ast z_1+{\dots}+z_d^\ast z_d$ of $\cA_0$ the algebraic set
$\cZ(f)$ corresponds to the spherical isometries in \cite{hmp}. Many
considerations on noncommutative real geometry based on {\it free}
polynomial algebras by J.W. Helton and his coworkers fit nicely into
the above concepts.

Let $f{=}(f_1,\dots,f_{d+1})$, where
$f_j(z_1,\dots,z_d){=}(1{-}z_j^\ast z_j)^2$ for $j{=}1,\dots,d$ and
$f_{d+1}(z_1,\dots,z_d){=}(1{-}\sum_{l=1}^d z_lz_l^\ast )^2$. Then
the elements of the algebraic set $\cZ(f)$ for $\cA_0$ is in
one-to-one correspondence with representations of the Cuntz algebra
$\cO_d$. \hfill $\circ$
\end{thex}

\subsection{The Role of Well-Behaved Unbounded Representations}

In this subsection we will show that in case of unbounded
$\ast$-representations one has to select "good"
$\ast$-representations rather than taking {\it all} (irreducible)
$\ast$-representations as $\cR$.

Let $\tau_{st}$ denote the finest
locally convex topology on a vector space.

\begin{thp}\label{closcone2}
If $\cA$ is the commutative $\ast$-algebra
$\dC[t_1,\dots,t_d]$, the Weyl algebra $\cW(d)$, the enveloping
algebra $\cE(\cG)$ or the free $\ast$-algebra $\dC\left<
t_1,\dots,t_d\right >$ (see Examples \ref{pol1}--\ref{free1}), then
the cone $\sum \cA^2$ is $\tau_{st}$-closed in $\cA$.
\end{thp}
\mn
{\bf Proof.} \cite{rev}, Theorem 4.2, p. 95, see e.g. \cite{unb}, Corollary 11.6.4. \hfill $\Box$
\begin{thp}\label{closcone3}
Let $\cA$ be a countably generated complex unital $\ast$-algebra such that
$\sum ~\cA^2$ is $\tau_{st}$-closed in $\cA$. For  any
$a\in \cA_h$ the following are equivalent:
\begin{itemize}
\item[\rm (i):] $a \in \sum ~\cA^2$.
\item[\rm (ii):] $ \pi(a)\geq 0$  for all $\ast$-representations $\pi$ of $\cA$.
\item[\rm (iii):] $ \pi(a)\geq 0$  for all irreducible $\ast$-representations $\pi$ of $\cA$.
\item[\rm (iv):] $f(a)\geq 0$ for each state $f$ of $\cA$.
\item[\rm (v):] $f(a)\geq 0$ for each pure state $f$ of $\cA$.
\end{itemize}
\end{thp}
\mn
{\bf Proof.} (i)$\to$(ii): $\langle \pi (\sum_j a_j^\ast a_j)\varphi,\varphi\rangle =
\sum_j \langle \pi(a_j)\varphi,\pi(a_j)\varphi \rangle \geq 0 $ for any $\varphi \in \cD(\pi)$.\\
(ii)$\to$(iv) and (iii)$\to$(iv): We apply (ii) resp. (iii) to the
GNS representation $\pi_f$ of the state $f$ and use formula (\ref{gns}). Note
that the GNS representation $\pi_f$ is irreducible if the state $f$
is pure (\cite{unb}, Corollary 8.6.8).\\
(iv)$\to$(i): Assume to the contrary that $a$ is not in $\sum
\cA^2$. By the separation theorem for convex sets (see e.g.
\cite{loc}, II.9.2), applied to the compact set $\{a\}$ and the
closed (!) convex set $\sum \cA^2$ of the locally convex space
$\cA_h[\tau_{st}]$, there exist a $\dR$-linear functional $g$ on
$\cA_h$ such that $g(a) <  {\rm inf}~ \{ g(c); c \in \sum \cA^2\}$.
Since $\sum \cA^2$ is a wedge, the infimum is zero, so we have
$g(a)<0$ and $g(\sum \cA^2)\geq 0$. The latter implies that the
Cauchy-Schwarz inequality holds. Therefore, $0\neq |g(a)|^2 \leq
g(1)g(a^2)$ which yields $g(1) >0$. Extending the $\dR$-linear
functional $g(1)^{-1}g$ on $\cA_h$ to a $\dC$-linear functional $f$
on $\cA$, we obtain a state $f$ such that $f(a) <0$.\\
(v)$\to$(iv): Since $\cA$ is countably generated, the assumptions of
Theorem 12.4.7 in \cite{unb} are satisfied. By this theorem, each
state  of $\cA$ is an integral over {\it pure} states.
This in turn gives the implication (v)$\to$(iv). \\
Since the implications (ii)$\to$(iii) and (iv)$\to$(v) are trivial,
the equivalence of (i)--(v) is proved. \hfill $\Box$

\mn Let $\cA$ be one of $\ast$-algebras from Proposition
\ref{closcone2}. Since then $\sum \cA^2$ is $\tau_{st}$-closed,
Proposition \ref{closcone3} applies and states that the sums of
squares in $\cA$ are precisely those elements which
 are nonnegative in {\it all} irreducible $\ast$-representations (or for
{\it all} pure states) of $\cA$. In particular, there is no
difference between the commutative $\ast$-algebra
$\dC[t_1,\dots,t_d]$ and the free $\ast$-algebra
$\dC\left<t_1,\dots,t_d\right>$ in this respect. In order to get
an interesting theory in the spirit of classical real algebraic
geometry one has to select a distinguished class $\cR$ of {\it
well-behaved} $\ast$-representations rather taking all irreducible
$\ast$-representations. For the $\ast$-algebras
$\dC[t_1,{\dots},t_d]$, $\cW(d)$ and $\cE(\cG)$ families $\cR$ of
such representations has been chosen in Examples
\ref{pol1}--\ref{env1}. It should be noted that there is no
general procedure for finding well-behaved representations of
arbitrary $\ast$-algebras.

Using essentially the $\tau_{st}$-closedness of the cone $\sum
\dC\left<t_1,{\dots},t_d\right>^2$ proved in \cite{rev}
we give a  short proof of the following theorem due to
Helton \cite{hel}.
\begin{thp}
Let $\cA{=}\dC\left< t_1,\dots,t_d \right>$ be the free complex
$\ast$-algebra in $d$ hermitian indeterminates $t_1,{\dots},t_d$
and
 $a \in \cA_h$. If $\pi(a) \geq 0$ for all finite dimensional
$\ast$-representations $\pi$ of $\cA$, then $a \in \sum \cA^2$.
\end{thp}
\mn {\bf Proof.} Let $\pi$ be a $\ast$-representation of $\cA$ and
$\varphi \in \cD(\pi)$. By Proposition \ref{closcone3},(ii)$\to$(i),
it suffices to show that there is a {\it finite dimensional}
$\ast$-representation $\rho$ such that $\varphi \in \cD(\rho)$ and
$\langle \pi(a)\varphi,\varphi \rangle = \langle \rho(a)
\varphi,\varphi \rangle$. This is easily done as follows.

Let $\cA_k$ be the vector space of polynomials
%in $t_1,{\dots},t_d$
of degree less than k. We choose $k$ such that $a \in \cA_k$.
Let $P$ denote the projection of the
%underlying
Hilbert space $\Hh$ on
the finite dimensional subspace $\pi(\cA_k)\varphi$. Since
$T_j:=P\pi(t_j)\lceil P\Hh$, $j{=}1,{\dots},d$, are selfajoint
operators on $P\Hh$, there is a $\ast$-representation $\rho$ of
$\dC\left<t_1,{\dots},t_d\right>$ on $P\Hh$ such
that $\rho(t_j)=T_j$. By construction we have $\pi(b)\varphi =\rho
(b) \varphi$ and hence $\langle \pi(b)\varphi,\varphi \rangle =
\langle \rho(b) \varphi,\varphi \rangle$ for all $b \in \cA_k$, so
in particular for $b=a$. \hfill $\Box$

\section{Positivstellens\"atze for General $\ast$-Algebras}

\subsection{Artin's theorem for General $\ast$-Algebras}\label{artin}

Let us begin our discussion with the commutative case. By Artin's
theorem for each nonnegative polynomial $a$ on $\dR^d$ there
exists a nonzero polynomial $c \in \dR[t]$ such that $ c^2 a
\in \sum \dR[t]^2$. For a noncommutative $\ast$-algebra $\cA$ a
natural guess is to generalize the latter to $c^\ast ac \in \sum
\cA^2$. (One might also think of $\sum_l c_l^\ast
ac_l \in \sum~\cA^2$, but Proposition \ref{ps3} below shows
that such a condition corresponds to a Nichtnegativstellensatz
rather than a Positivstellensatz.)

In the
commutative case the relation $c^2 a \in \sum \dR[t]^2$ implies that
the polynomial $a$ is nonnegative on $\dR^d$. However, in the
noncommutative case such a converse is not true in general as the
following  examples show.

\begin{thex}
Let $\cA$ be the Weyl algebra $\cW(1)$ and $\cR{=}\{\pi_0 \}$, see
Example \ref{weyl1}. Set $N{=}a^\ast a$. Since  $aa^\ast{-}a^\ast
a{=}1$,  we have $a(N{-}1)a^\ast {=} N^2{+} a^\ast a \in \sum
\cA^2$. But $\pi_0 (N{-}1)$ is not nonnegative, since $\langle
\pi_0(N{-}1)e_0,e_0\rangle =-1$ for the vacuum vector $e_0$. \hfill
$\circ$
\end{thex}

\begin{thex}
Let $\cA$ be the $\ast$-algebra with a single generator $a$ and defining relation $a^\ast
a=1$. Then $p_0:=1{-}aa^\ast$ is a nonzero projection in $\cA$ and we have
$p_0ax a^\ast p_0 =
 0 \in \sum \cA^2$ for arbitrary
$x\in \cA$. But elements of the form $axa^\ast$ are in general not
nonnegative in $\ast$-representations of $\cA$. \hfill $\circ$
\end{thex}

For a reasonable generalization of Artin's theorem one should add
conditions which ensure that $\pi(a) \geq 0$ for $\pi \in \cR$. In
the commutative case $c$ can be chosen such that the zero set
$\cN(c)$ is contained in the zero set $\cN (a)$. (This follows, for
instance, from Stengle's Positivstellensatz.)  It seems to be
natural to require a generalization of this condition in the
noncommutative case as well. Thus, our {\it first version of a
noncommutative generalization of Artin's theorem} for $\cA$ and
$\cR$ is the following assertion:

{\it For each   $a \in \cA_h$ such that $\pi (a) \geq 0$ for all
$\pi \in \cR$  there exists an element  $c\in \cA$ such that}
\begin{align}\label{suma}
&c^\ast ac \in \sum \cA^2,\\
&\cN(\pi(c)^\ast) \subseteq \cN(\overline{\pi (a)})~~{\rm for~~each} ~~ \pi \in \cR.\label{kerc}
\end{align}
Let $a\in \cA_h$ and suppose conversely that there exists a $c \in
\cA$ such that (\ref{suma}) and (\ref{kerc}) hold. For $\pi \in \cR$
we put $\cE_\pi:= \pi(c)\cD(\pi) + \cN(\overline{\pi (a)})$.
\begin{thl}
$\cE_\pi$ is dense in $\Hh(\pi)$ and
$\langle \ov{\pi (a)}\eta, \eta \rangle \geq 0$ for $\eta \in \cE_\pi$.
\end{thl}
\mn {\bf Proof.} Since $\Hh = \cR(\ov{\pi(c)})\oplus
\cN(\pi(c)^\ast) $ and $\cR(\pi(c))$ is dense in
$\cR(\ov{\pi(c)})$,the linear subspace $\cE_\pi$ is dense in
$\Hh(\pi)$. For $\varphi \in \cN(\pi (c)^\ast)$ and $\psi \in
\cD(\pi)$, using condition (\ref{kerc}) we obtain
$$
\langle \ov{\pi(a)}(\varphi + \pi(c)\psi), \varphi + \pi(c)\psi \rangle =
\langle \pi(a)\pi(c)\psi, \pi(c)\psi \rangle = \langle \pi(c^\ast ac)\psi, \psi \rangle \geq 0,
$$
where the last inequality follows at once from condition (\ref{suma}). \hfill $\Box$

\mn
Since $\cE_\pi$ is dense in $\Hh(\pi)$, it is obvious that $\pi(a)
\geq 0$ on $\cD(\pi)$ when the operator $\pi(a)$ is {\it bounded}.
If $\pi(a)$ is unbounded, it follows that $\pi(a)\geq 0$ on
$\cD(\pi)$ if we replace (\ref{kerc}) by the following
technical condition:
\begin{align}\label{kerc^}
\cN(\ov{\pi(a)})+ \pi(c)\cD(\pi)~~{\rm is~ a ~core~ for}~ \ov{\pi(a)}.
\end{align}
In many cases it is difficult to decide whether or nor (\ref{kerc})
can be satisfied. We now formulate another  condition which is often
easier to verify.

Let $\Psi(c)$ denote the set of all finite sums of linear
functionals of the form $\langle
\pi(\cdot)\pi(c)\varphi,\pi(c)\varphi \rangle$ on $\cA$, where $\pi
\in \cR$ and $\varphi \in \cD(\pi)$. By our {\it second version of a
noncommutative generalization of Artin's theorem} we mean that
assertion (\ref{suma}) and the following density condition
(\ref{weakden}) hold:
\begin{align}\label{weakden}\nonumber
&{\rm For~each}~ \pi \in \cR~{\rm and}~ \psi \in \cD(\pi) ~{\rm the
~functional}~ \langle \pi(\cdot)\psi,\psi \rangle \\ ~&{\rm
on}~\cA~{\rm is~the~weak~ limit~ of~a~net~of~ functionals~ from}~
\Psi(c).
\end{align}
Since (\ref{suma}) obviously implies that $\pi(a) \geq 0$ on
$\pi(c)\cD(\pi)$, it follows from (\ref{suma}) and (\ref{weakden})
that $\pi(a) \geq 0$ on $\cD(\pi)$ for all $\pi \in \cR$. Clearly,
if suffices to assume (\ref{weakden}) for vectors $\psi$ from a core
for $\pi(a)$ rather than for all $\psi \in \cD(\pi)$.

We give an example for this second version. This example is due to
Y. Savchuk and details of proofs will appear in his forthcoming
thesis.
\begin{thex}\label{yurii}
Let $\cA$ be the complex $\ast$-algebra with a generator
$a$ and defining relation $a^\ast a+aa^\ast =1$. All $\ast$-representations of $\cA$ act
by bounded operators. Let $\cR$ be the equivalence classes of all
irreducible $\ast$-representations. They are formed by series
$\rho_{\alpha,\varphi}$, where $\alpha \in [0,1/2), \varphi \in
[0,2\pi)$, of 2-dimensional representations and $\rho_\varphi$,
where $\varphi \in [0,2\pi)$, of 1-dimensional representations.
These representations act on the generator $a$ by
$$\rho_{\alpha,\varphi}(a)=
\left(
  \begin{array}{cc}
    0 & e^{\ii\varphi}\sqrt{\alpha} \\
    \sqrt{1-\alpha} & 0 \\
  \end{array}
\right)~~~ \mbox{and}~~~~ \rho_\varphi(a) =\frac{e^{\ii\varphi}}{\sqrt{2}}.
$$
For $\alpha{=}1/2$ the matrix $\rho_{\alpha,\varphi}(a)$
defines {\it reducible} ´$\ast$-representation of
$\cA$.
For the $\ast$-algebra $\cA$ we have the following Positivstellensatz:\\
{\it Suppose that $b\in \cA_h$ and $\pi(b)\geq 0$ for all $\pi \in
\cR$. Then there exists an element $c=c^\ast$ of the center of $\cA$ such
that $c^2 a \in \sum ~\cA^2$ and condition (\ref{weakden}) is satisfied.}

Let $\cB$ denote the $\ast$-algebra of  complex polynomials in
three commuting indeterminates $u{=}u^\ast, v, v^\ast$ satisfying
the relation $u^2+vv^\ast =1$. The map
$$a \to
\left(
  \begin{array}{cc}
    0 & v \\
    u & 0 \\
  \end{array}
\right)
$$
extends to a $\ast$-isomorphism of $\cA$ onto a $\ast$-subalgebra of
the matrix algebra $\cM_2(\cB)$. If we consider $\cA$ as a
$\ast$-subalgebra of $\cM_2(\cB)$, then the element $c$ is a
multiple $c_0{\cdot}I$, where $c_0 \in \cB_h$, of the unit matrix
$I$. \hfill $\circ$
\end{thex}

\subsection{Generalizations of Stengle's Theorem to General $\ast$-Algebras}\label{st}

As already noted in subsection \ref{ncsas} the usual definition of
the preorder does not make sense in the noncommutative case, because
the product of noncommuting hermitian elements is not hermitian. For
arbitrary $\ast$-algebras and semialgebraic sets I don't know how a
proper generalization of the preorder might look like. In this
subsection we propose one possible way to remedy this difficulty by
reducing the problem to some appropriate {\it commutative}
$\ast$-subalgebra. Our guiding examples for this method are
$\ast$-subalgebras of matrix algebras $\cM_n(\cB)$ (Example
\ref{matcomm}).

Let $\cZ(\cA)$ be the center of $\cA$ and let $\sum_{\cZ}$ denote the set of nonzero elements of the wedge
$\sum \cZ(\cA)^2$. We shall assume the
following:\\
 {\it If $az=0$ for some $a \in \cA$ and  $z \in \cZ(\cA)$,
then $a=0$ or $z=0$.}

Obviously, this  is fulfilled if $\cA$ has no
zero divisors.
%and it implies that the center $\cZ(\cA)$ has no zero divisors.
\begin{thd}\label{def1}
For $a_1,a_2 \in \cA_h$, we  write $a_1 {\sim} a_2$ if there
exist elements $s_1,s_2 \in \sum_{\cZ}$, $z \in \cZ(\cA)$ and $x_\pm \in \cA$ such that
\begin{align}\label{equiv}
x_- x_+ =x_+x_-=z~~{\rm and}~~
 s_1~ a_1 =s_2 ~x_+a_2x_+^\ast.
\end{align}
\end{thd}

\begin{thl}\label{equivlemma}
$"\sim"$ is an equivalence relation on $\cA_h$.
\end{thl}
\mn {\bf Proof.} Suppose  $a_1 \sim a_2$. Multiplying the second
equation of (\ref{equiv}) by $x_-$ from the left and by $x^\ast_-$
from the right and using the first equations  we get
$$
(s_2 zz^\ast)~ a_2 =s_1~x_-a_1x_-^\ast.
$$
Since $s_2zz^\ast \in \sum_{\cZ}$, the latter means that $a_2 \sim a_1$.

Suppose  $a_1 \sim a_2$ and $a_2 \sim a_3$. Then there are elements
$s_1,s_2, s_3, s_4 \in \sum_{\cZ}$, $x_+,x_-, y_+, y_- \in \cA$ and
$z,w \in \cZ(\cA)$ such that $x_-x_+=x_+x_-=z$, $y_-y_+=y_+y_-=w$,
$s_1~ a_1 =s_2 ~ x_+a_2 x_+^\ast$, $s_3~a_2= s_4~ y_+a_3 y_+^\ast$.
Setting $u_-:=y_-x_-$, $u_+:=x_+y_+$, we have $u_-u_+= u_+u_-=zw$
and $s_3s_1 ~a_1 =s_2x_+(s_3a_2)x_+^\ast = s_2s_4 x_+y_+ a_3
y_+^\ast x_+^\ast = s_2s_4 u_+ a_3 u_+^\ast$, so that $a_1 \sim
a_3$.

Since obviously $a \sim a$,  $"\sim"$ is an equivalence relation.
\hfill $\Box$

\mn Let $\cC$ be a quadratic module of $\cA$ and assume that $a_1
\sim a_2$. Then we have $s_1a_1 \in \cC$ if and only if $s_2zz^\ast
a_2 \in \cC$. That is, up to multiples from the set $\sum_{\cZ}$,
$a_1$ belongs to $\cC$ if and only if $a_2$ is in $\cC$.

The relation $\sim$ is extended to tuples  $a{=}(a_1,\dots,a_n)$ and
$b{=}(b_1,\dots,b_r)$ from $\cA_h$ by defining  $a
\sim b$  if $a_j \sim b_l$  for all $j{=}1,{\dots},n$ and $
l{=}1,{\dots},r$.
\begin{thd}\label{def2}
Suppose  $a \sim b$. We shall write $a \sim^+ b$ if for any
representation $\pi \in \cR$, $\pi (a_j) \geq 0$  for all
$j{=}1,{\dots},n$ implies that $\pi (b_l)\geq 0$ for all
$l{=}1,{\dots},r$ and we write $a \overset{+}{\sim} b$ if $a \sim^+ b$  and
$b \sim^+ a$.
\end{thd}

To begin with the setup for Stengle's theorem, let us fix a k-tuple
$f=(f_1,\dots,f_k)$ of  elements $f_j \in \cA_h$. Let $a$ be an
element of $\cA_h$ which is nonnegative on the semialgebraic set
$\cK(f)$, that is, $\pi (a) \geq 0$ for $\pi \in \cK(f)$.

Suppose there exist a finitely generated {\it commutative} real
subalgebra $\cB$ of $\cA_h$ such that the following assumptions are fulfilled:

\mn
(I)~~{\it There exist a finite tuples $c=(c_1,\dots,c_m)$ and
$b=(b_1,\dots,b_r)$ of elements of $\cB$ such that $a
\overset{+}{\sim} c$ and $f \overset{+}{\sim} b$.}\\
(II) ~{\it For $j{=}1,{\dots},m$, we have $\pi (c_j)\geq 0$ for all $\pi \in
\cK(b)$ if and only if $c_j(s)\geq 0$ for all $s \in \cK_b$.}

\mn Recall that $\cK(b)$ is the noncommutative semialgebraic set
defined by (\ref{bas}) and $\cK_b = \{ s \in \hat{\cB}: b_1(s)\geq
0,{\dots}, b_r(s)\geq 0 \}$ is the "ordinary" semialgebraic set for
the commutative real algebra $\cB$ defined by (\ref{kf}).

We now derive our noncommutative version of Stengle's theorem. By
assumption (I), we have $f \overset{+}{\sim} b$ and $a \sim^+ c$.
The relation  $f \overset{+}{\sim} b$ implies that the two
semialgebraic sets $\cK(f)$ and $ \cK(b)$ of $\cA$ coincide. Since
$a \sim^+ c$, we have $\pi (c_j)\geq 0$ for all $\pi \in
\cK(f)=\cK(b)$ and $j{=}1,{\dots},m$. Therefore, by assumption
(II), $c_j \geq 0$ on $\cK_b$. Let $\cT_b$ denote the preorder
(\ref{tf}) for the commutative algebra $\cB$. By Stengle's
theorem, applied to $\cK_b$ and $\cT_b$, there exist elements
$g_j, h_j \in \cT_b$ and numbers $n_j \in \dN$ such that
\begin{align}\label{stenglea1}
g_j c_j = c_jg_j =c^{2n_j}_j + h_j.
\end{align}
Since $c_j \sim a$ and $b_l \sim f_l$, there exist elements
$s_{1j}, s_{2j}, s_{3l}, s_{4l} \in \sum_{\cZ}$ and $x_{+j},
y_{+l}\in \cA$ such that $s_{1j} c_j = s_{2j} x_{+j}a x_{+j}^\ast$
and $s_{3l}b_j = s_{4l} y_{+l}f_ly_{+l}^\ast$. Put $s_3
{:=}s_{31}\cdots s_{3r}$. Multiplying (\ref{stenglea1}) by the
central element $s_3s_{1j}^{2n_j+1}$ we obtain
\begin{align}\label{stenglea2}
s_3 g_j~ s_{1j}^{2n_j} s_{2j} ~ x_{+j} ax_{+j}^\ast = s_3
s_{1j}(s_{2j} x_{+j}ax_{+j}^\ast)^{2n_j} + s_{1j}^{2n_j+1}~ s_3h_j.
 %j=1,{\dots},m.
\end{align}
for $j{=}1,{\dots},m$. Set $p_j{:=}s_3 g_j s_{1j}^{2n_j}s_{2j}$. Let
$\cT(f)$ denote the quadratic module of $\cA$ generated by the
preorder $\cT_b$ of $\cB$. Since
 $s_3 g_j$ and $s_3 h_j$ belong to $\cT(f)$, $p_j
\in \cT(f)$. Hence the right-hand side of (\ref{stenglea2}) is in
$\cT(f)$, so we have
\begin{align}\label{stenglea3}
p_j x_{+j}a x_{+j}^\ast =x_{+j}a x_{+j}^\ast p_j  \in \cT(f) ~{\rm for} ~j=1,\dots,m.
\end{align}
That is, there exist elements $p_j \in \cT(f)$ and $x_{+j} \in \cA$
such that (\ref{stenglea3}) holds. We consider this statement (and
likewise the more precise  equalities (\ref{stenglea2}))  as a {\it
noncommutative version of Stengle's theorem}.

We now turn to the converse direction, that is, we show that our version of
Stengle's theorem implies that $\pi(a)\geq 0$ for all $\pi \in \cK(f)$
. For suppose that (\ref{stenglea2}) is satisfied for
$j=1,{\dots},m$ with $s_{1j},s_{2j},s_3 ,g_j,h_j ,x_{+j}$ as above.
Then, since (\ref{stenglea2}) is nothing but (\ref{stenglea1})
multiplied by $s_{1j}^{2n_j+1} s_3$, equation (\ref{stenglea1})
holds, so  each element $c_j$ is nonnegative on the set $\cK_b$.
Therefore, by assumption (II) we have $\pi (c_j)\geq 0$ for all $\pi
\in \cK(b)=\cK(f)$. Since $c \sim^+ a$ by assumption (I),  it
follows that $\pi (a) \geq 0$ for all $\pi \in \cK(f)$.

\mn We close this subsection by discussing assumption (II). The
following simple example shows that it is not always satisfied.

\begin{thex}
Let $\cB$ be the $\ast$- algebra $\dC[t]$ of complex polynomials in
one hermitean indeterminate $t$. From the moment problem theory it
is known that there exists a state $f$ on $\cB$ such that $f(t^3p
\bar{p})\geq 0$ for all $p \in \cB$ and $f(t p_0^2)<0$ for some $p_0
\in \cB_h$. For the GNS representation $\pi_f$ of $f$ we then have
$\pi_f(t^3)\geq 0$ and $\pi_f(t)\not\ge 0$. Therefore, if the family
$\cR$ contains $\pi_f$, then $t \not\in \cK(t^3)$, but $t \in
\cK_{t^3}=[0,\infty)$.

The converse direction fails if $\cR$ is to "small". For instance,
if we take $\cR{=} \{\pi_s; s\in [0,1] \}$, where $\pi_s(p)=p(s)$,
then $t{-}1 \in \cK(t)$, but $t{-}1 \not\in \cK_t$. \hfill $\circ$
\end{thex}

To give a sufficient condition for assumption (II), we fix
hermitian generators $y_1,\dots,y_d$ of $\cB$. Then
 $\hat{\cB}$ becomes a subset of $\dR^d$ by
identifying a character with its values at the generators. If $\pi$
is a bounded $\ast$-representation of $\cB$ on a Hilbert space
, then the $d$-tuple of commuting bounded selfadjoint operators
$\pi(y_j)$ has a unique spectral measure $E_\pi$. \mn
\begin{thp}
Suppose  $\cR$ is family of bounded $\ast$-representation of
$\cB$ on Hilbert spaces. Let $b{=}(b_1,\dots,b_r)$ be an $r$--tuple of
elements of $\cB$ such that $\cK_b$ is the union of all its subsets of
the form ${\rm supp}~E_\pi$, where $\pi \in \cR$. Then for any $a
\in \cB$, we have $a \in \cK(b)$ if and only if $a \in \cK_b$.
\end{thp}
\mn {\bf Proof.} Since $\pi \in \cR$ is bounded, it is a direct
sum of cyclic representations. Hence we can assume without loss of
generality that each $\pi \in \cR$ has a cyclic vector
$\varphi_\pi$. Then $\mu_\pi(\cdot):=\langle E_\pi(\cdot)
\varphi_\pi , \varphi_\pi \rangle$ defines a positive Borel
measure on $\dR^d$ such that ${\rm supp}~ E_\pi = {\rm supp}~
\mu_\pi$. From the spectral theorem we obtain
%for any $p \in \dC[y_1,\dots,y_d]$ we have
\begin{align}\label{int}
\langle\pi(p(y))\varphi_\pi ,\varphi_\pi \rangle = \tint_{\dR^d}~ p(s)~ {\rm d}\mu_\pi (s)~
{\rm for}~ p\in \dC[y_1,...,y_d].
\end{align}
Suppose that $c \in \cB_h\cong \dR[y_1,\dots,y_d]$. For $p \in \dC[y_1,\dots,y_d]$, we have
$$
\langle \pi(c)\pi(p)\varphi_\pi,\pi(p) \varphi_\pi \rangle =
\tint_{\dR^d}~ c(s) |p(s)|^2  ~{\rm d}\mu_\pi (s).
$$
Since the polynomials are uniformly dense in the continuous function
on the compact set $supp \mu_\pi$, it follows that $\pi (c) \geq 0$
if and only if ${\rm supp}~ \mu_\pi \subseteq \cK_c$. This implies
that $\cK(b)= \{\pi \in \cR: {\rm supp}~ \mu_\pi \subseteq \cK_b
\}$. Therefore, if $a \in \cK_b$, then $a \in \cK(b)$ by
(\ref{int}). Conversely, if $a \in \cK(b)$, then ${\rm supp}~\mu_\pi
\subseteq \cK_a$ for all $\pi \in \cR$ such that ${\rm supp}~\mu_\pi
\subseteq \cK_b$. By assumption $\cK_b$ is the union of all sets
${\rm supp}~\mu_\pi={\rm supp}~ E_\pi $ which are contained in
$\cK_b$. Hence $\cK_b \subseteq \cK_a$ which in turn yields $a \in
\cK_b$. \hfill $\Box$

\subsection{Diagonalization of Matrices with Polynomial Entries}
In the rest of this section, $\cA$ is the real $\ast$-algebra
$\cM_n(\dR[t])$ of $n\times n$-matrices over
$\dR[t]=\dR[t_1,\dots,t_d]$ with involution given by the transposed
matrix $A^t$ of $A$ and $\cR$ is the set $\{\rho_s; s \in \dR^d\}$
of irreducible $\ast$-representations $\rho_s{:}A \to A(s)$, see
Example \ref{matcomm}. Then $\cA_h$ is the set $\cS_n(\dR[t])$ of
symmetric matrices and the unit of $\cA$ is the unit matrix $I$.
Clearly, $\rho_s(A)=A(s)\geq 0$ if and only if the matrix $A(s)$ is
positive semidefinite.

We begin with some notation. If $i{=}(i_1,{\dots},i_p)$ and
$j{=}(j_1,{\dots},j_p)$ are $p$-tuples of integers such that $1\leq
i_1< \dots < i_p\leq n$ and $1\leq j_1 < \dots < j_p\leq n$ and $A
\in \cA$, then $M^i_j{=}M^i_j(A)$ denotes the principal minor of $A$
with columns $i_k$ and rows $j_k$. If $i_1{=}j_1{=}1,\dots,,
i_p{=}j_p{=}p$, we write $M_p$ instead of $M_j^i$.

For $\lambda =(\lambda_1,\dots,\lambda_p)$, where
$p \le n$, let $D(\lambda)=D(\lambda_1,\dots,\lambda_p)$ denote
the $n\times n$ diagonal matrix  with diagonal entries
$\lambda_1,\dots,\lambda_p,0,\dots,0$.

Now let $A \in \cS_n(\dR[t])$, $A\neq 0$, $n\geq2$, and assume that $A$ has
rank $r$ and that $M_1(A)\neq 0,\dots, M_r(A)\neq 0$. If the
latter is true, we say that $A$ has {\it standard form}. For such
a matrix $A$ we define two lower triangular $n\times n$-matrices
$Y_{\pm} =(y_{ij}^{\pm})$ with entries given by the rational
functions
\begin{align*}
&y_{ij}^\pm =\pm M^{(1,\dots,j-1,i)}_{(1,\dots,j-1,j)} M_j^{-1}~
{\rm for}~ j{=}1,\dots,r,~i{=}j+1,\dots,n,\\
&y_{ii}^\pm =1~{\rm for}~ i{=}1,\dots,n,\\
&y_{j}^\pm =0 ~{\rm otherwise}~(j{=}r+1,\dots,n,i{=}j+1,\dots,n~{\rm and}~ i\geq j,i,j{=}1,\dots,n).
\end{align*}
By Satz 6.2 in \cite{gam}, p.64, we have
$A=Y_+D(M_1,M_2M_1^{-1},\dots,M_rM_{r{-}1}^{-1}) Y_+^t$. Since
obviously $Y_+^{-1}=Y_-$, the latter yields
\begin{align}\label{DYA}
D(M_1,M_2M_1^{-1},\dots,M_rM_{r-1}^{-1}) = Y_-AY_-^t.
\end{align}
Set $D{=}M_1{\cdots} M_{r{-}1} D(M_1,M_2M_1^{-1},{\dots},M_rM_{r{-}1}^{-1})$ and
$X_\pm {=}M_1{\cdots} M_{r{-}1} Y_\pm$. Clearly, $D$ and $X_\pm$ are
in $\cM_n(\dR[t])$. From the relations $Y_-^{-1}=Y_+$ and
(\ref{DYA}) we obtain
\begin{align}\label{X}
X_+X_- &=X_-X_+ =(M_1\cdots M_{r{-}1})^2 I,\\
(M_1\cdots M_{r{-}1})^4 A&= X_+DX_+^t,~~ D= X_- AX_-^t.\label{DA}
\end{align}
That is, we have
shown that {\it for any matrix $A \in \cS_n(\dR[t])$ in standard
form (that is, rank $A =r$ and $M_1\neq 0,\dots,M_r\neq 0$) there
exists a diagonal matrix $D\in \cM_n(\dR[t])$ such that $A \sim D$
and the corresponding matrices $X_\pm$ can be chosen to be lower
triangular.}

We now turn to arbitrary matrices in $\cS_n(\dR[t]).$ Our aim is
to prove Proposition \ref{diagmat} below. The main technical
ingredient for this proof is the following procedure for block
matrices over a ring $R.$ We write a matrix $A\in\cS_n(R),\ n\geq
2,$ as
$$A=
\left(
  \begin{array}{cc}
    \alpha & \beta \\
    \beta^t & C \\
  \end{array}
\right), \mbox{where}\ C\in\cS_{n-1}(R),\ \beta\in\cM_{1,n-1}(R),
$$
and put
$$X_{\pm}=
\left(
  \begin{array}{cc}
    \alpha & 0 \\
    \pm\beta^t & \alpha I \\
  \end{array}
\right),\ \widetilde{A}=\left(
                           \begin{array}{cc}
                             \alpha^3 & 0 \\
                             0 & \alpha(\alpha C-\beta^t\beta) \\
                           \end{array}
                         \right).
$$
Then we have
\begin{align}\label{vor_diagmat_1}
X_+X_-=X_-X_+=\alpha^2\cdot I,\\
\alpha^4 A=X_+\widetilde{A}X_+^t,\
\widetilde{A}=X_-AX_-^t.\label{vor_diagmat_2}
\end{align}

\begin{thp}\label{diagmat}
Let $A \in \cS_n(\dR[t])$,$A\neq 0$. Then there exist diagonal
matrices $D_l \in \cM_n(\dR[t])$, matrices $X_{\pm,l} \in
\cM_n(\dR[t])$ and polynomials $z_l \in \sum \dR[t]^2$,
$l=1,\dots,m,$ such that:
\begin{itemize}
\item[\em (i)] $X_{+l}X_{-l} =X_{-l} X_{+l} =z_j I,~ D_l =X_{-l}AX_{-l}^t, ~
z_l A =X_{+l}D_l X_{+l}^t$,
\item[\em (ii)] For $s \in \dR^d$,~ $A(s) \geq 0$ if and only if~ $D_l(s)\geq 0$ for all
$l=1,\dots,m$.
\end{itemize}
\end{thp}%%%%%%%%%%%%%%%%%%%%%%%%%%% EINSCHUB 2
\mn {\bf Proof.} Let $i,j\in\{1,\dots,n\},\ i\leq j.$ Put
$\widetilde{a}_{ii}=a_{ii}$ and
$\widetilde{a}_{ij}=a_{ij}+\frac{1}2(a_{ii}+a_{jj})$ if $i<j.$ We
first show that there is an orthogonal matrix $T_{ij}\in\cM_n(\dR)$
such that
\begin{align}\label{diagmat_proof_1}
T_{ij}AT_{ij}^t=\left(
                  \begin{array}{cc}
                    \widetilde{a}_{ij} & \ast \\
                    \ast & \ast \\
                  \end{array}
                \right).
\end{align}

For $l\in\{2,\dots,n\},$ let $P_l$ denote the permutation matrix
which permutes the first row and the $l$-th row. Setting $T_{11}=I$
and $T_{ii}=P_i$ for $i=2,\dots,n$, (\ref{diagmat_proof_1}) holds for
$i=j.$ Now suppose $i<j.$ Let $S=(s_{kl})\in\cM_n(\dR)$ be the
matrix with $s_{ii}=s_{ij}=s_{ji}=2^{-1/2},\ s_{jj}=-2^{1/2},\
s_{ll}=1$ if $l\neq i,j$ and $s_{kl}=0$ otherwise. Set
$T_{ij}=T_{ii}S.$ One easily checks that $T_{ij}$ is orthogonal and
(\ref{diagmat_proof_1}) is satisfied for $i<j.$

Now we apply the above procedure to the block matrix
$A_{ij}:=T_{ij}AT_{ij}^t$, $ i\leq j,$ given by
(\ref{diagmat_proof_1}). Let $\widetilde{A}_{ij},\ X_{\pm,ij}$ denote
the corresponding matrices. Then there is a matrix
$B_{ij}\in\cS_{n-1}(\dR[t])$ such that
\begin{align}\label{diagmat_proof_2}
\widetilde{A}_{ij}=\left(
                     \begin{array}{cc}
                       \widetilde{a}_{ij}^3 & 0 \\
                       0 & B_{ij} \\
                     \end{array}
                   \right).
\end{align}
We claim that for any $s\in\dR^d,\ A(s)\geq 0$ if and only if
$\widetilde{a}_{ij}(s)\geq 0$ and $B_{ij}(s)\geq 0$ for all
$i,j\in\{1,\dots,n\},\ i\leq j.$

Indeed, if $A(s)\geq 0,$ then $A_{ij}(s)\geq 0$ by
(\ref{diagmat_proof_1}) and hence
$\widetilde{A}_{ij}=X_{-,ij}A_{ij}(s)X_{-,ij}^t\geq 0$ by
(\ref{vor_diagmat_2}), so $\widetilde{a}_{ij}(s)\geq 0$ and
$B_{ij}(s) \geq 0$ by (\ref{diagmat_proof_2}). Conversely, assume
that $\widetilde{a}_{ij}(s)\geq 0$ and $B_{ij}(s) \geq 0$ for all
$i,j$, $i\leq j$. Then $\widetilde{A}_{ij}(s)\geq 0$ for all $i,j$.
If $\widetilde{a}_{ij}(s)=0$ for all $i,j$, then $a_{ij}(s)=0$ for
all $i,j$ and hence $A(s)=0$. If $\widetilde{a}_{ij}(s)>0 $ for some
$i,j$, we conclude that
$A_{ij}(s)=\widetilde{a}_{ij}(s)^{-4}X_{+,ij}\widetilde{A}_{ij}(s)X_{+,ij}^t
\geq 0$ by (\ref{vor_diagmat_2}) and so
$A(s)=T_{ij}^tA_{ij}(s)T_{ij}\geq 0$. This completes the proof of
the claim.

Applying the same reasoning to the matrices $B_{ij}$ instead of $A$
and proceeding by induction we obtain after at most $n{-}1$ steps a finite
sequence of diagonal matrices having the desired properties.\hfill
$\Box$
\begin{thc}\label{diagmat2}
For each matrix $A\in\cS_n(\dR[t]), A{\neq} 0$, there exist
nonzero polynomials $ b, d_j \in \dR[t]$, $j{=}1,\dots,r,r\leq n,$ and matrices
$X_+,\ X_- {\in} \cM_n(\dR[t])$ such that
$$
X_+X_-=X_-X_+=bI,\ b^2A=X_+DX_+^t,\ D=X_-AX_-^t,
$$
where $D$ is the diagonal matrix $D(d_1,\dots,d_r).$ In particular, $A\sim
D.$
\end{thc}
\mn {\bf Proof.} Since $A\neq 0,\ a_{ij}\neq 0$ for some $i,j$. We
apply the above procedure to the matrix $B_{ij}$ from
(\ref{diagmat_proof_2}) and proceed by induction until the
corresponding matrix $B_{ij}$ is identically zero.
\hfill $\Box$

\mn Remark.  Suppose that $A\in\cS_n(\dR[t]),\ A\neq 0.$ Because the
rank $r$ of $A$ is the column rank and the row rank it follows that
$A$ has a non-zero principal minor of order $r.$ Hence, there exists
a permutation matrix $P$ such that $M_r(PAP^t)\neq 0.$ But it may
happen that all principal minors of order $r{-}1$ of $PAP^t$ vanish,
so $PAP^t$ is \textit{not} in standard form. A simple example is
$$
\left(
  \begin{array}{ccc}
    1 & -1 & 1 \\
    -1 & 1 & 1 \\
    1 & 1 & 1 \\
  \end{array}
\right).
$$

\subsection{Artin's Theorem and Stengle's Theorem for Matrices of Polynomials}\label{artst}

From Corollary \ref{diagmat2} and Proposition \ref{diagmat} we
easily derive versions of Artin's theorem and Stengle's theorem
for matrices of polynomials.

The next proposition is Artin's theorem for matrices of
polynomials. It was first proved in \cite{gr} and somewhat later also
in \cite{ps}.
\begin{thp}\label{artin1}
Let $A \in \cS_n(\dR[t])$. If $A(t) \geq 0$ for all $t \in \dR^d$,
then there exist a polynomial $c\in \dR[t]$, $c{\neq} 0$, such
that $c^2 A \in \sum ~ \dR[t]^2$.
\end{thp}
\mn {\bf Proof.} Let $D=D(d_1,\dots,d_r)$ be the diagonal matrix
from Corollary \ref{diagmat2}. Since $D=X_-A X^t_-$ and $A(t)\geq
0$ on $\dR^d$, we have $D(t)\geq 0$ and hence $d_j(t) \geq 0$ on
$\dR^d$. Since $b^2 A=X_+DX_+^t$, the assertion follows at once by
applying Artin's theorem for polynomials to the diagonal entries
$d_1,\dots,d_r$ and multiplying by the product of denumerators.
\hfill $\Box$

\mn Let $c \in \dR[t]$, $c \neq 0$. Since the set $\{s\in
\dR^d:c(s)\neq 0\}$ is dense in $\dR^d$, each $s \in \dR^d$ is limit
of a sequence of points $s_n$ such that $c(s_n) \neq 0$. Then each
vector state of $\rho_s$ is weak limit of vector states of
$\rho_{s_n}$ with vectors from $c(s_n)\dR^d$. Since these
functionals belong to the set $\Psi(c)$,  condition (\ref{weakden})
is fulfilled and the second version of Artin's theorem  holds.

We now turn to Stengle's theorem and apply the setup of subsection
\ref{st} to the $\ast$-algebra $\cA=\cM_n(\dR[t])$ and its
commutative $\ast$-subalgebra $\cB$ of diagonal matrices.  Let
$F=(F_1,\dots,F_k)$ be a k-tuple of elements from $\cS_n(\dR[t])$
and let $\cK(F){=}\{\rho_s: s \in \dR^d, F_1(s)\geq 0,\dots, F_k(s)
\geq 0\}$ be the corresponding noncommutative semialgebraic set.
Suppose that $A \in \cS_n(\dR[t])$ and $\rho_s(A)=A(s) \geq 0$ for
all $s \in \cK(F)$.

By Proposition \ref{diagmat} there exists an $m$-tuple
$C=(C_1,\dots,C_m)$ of diagonal matrices such that $A
\overset{+}{\sim} C$. Applying Proposition \ref{diagmat} to each
matrix $F_j$ we obtain a finite sequence of diagonal matrices. Let
$B=(B_1,\cdots,B_r)$ denote r-tuple formed by all these diagonal
matrices and all diagonal matrices obtained by permutations of their
diagonal entries for $j{=}1,\dots,k$. By Proposition \ref{diagmat},
we then have $F \overset{+}{\sim}B$, so assumption (I) is satisfied.
The set $\hat{\cB}$ of characters of $\cB$ consists of all
functionals $h_{i,s}$, where $s\in \dR^d$ and $j{=}1,\dots,n$, given
by $h_{j,s}(D(d_1,\dots,d_n))= d_j(s)$. Let $b_{j1},\dots,b_{jn}$ be
the diagonal entries of $B_j$. Since $B$ contains all permuted
diagonal matrices, the set $\cK{:=}\{s\in \dR^d: b_{1l}(s)\geq
0,\dots,b_{rl}(s)\geq 0\}$ does not depend on $l{=}1,\dots,n$. Hence
we have $\cK(B)=\{ \rho_s: s \in \cK\}$ and $\cK_B=\{h_{i,s}:s \in
\cK \}$ which  implies that assumption (II) is fulfilled. Therefore
the version of Stengle's theorem stated in subsection \ref{st} is
valid. Recall that $\cT(F)$ is the quadratic module of $\cA$
generated by all products $B_{i_1}\cdots B_{i_l}$, where $1\leq
i_1<i_2<\dots i_l\leq r$. We may consider $\cT(F)$ as a
noncommutative substitute of the preorder associated with $F$. Since
$\cT(F)$ depends on the particular diagonalizations of $F_j$, it is
neither uniquely nor canonically associated with $F$.

\section{Archimedean Quadratic Modules}
In this section $\cA$ is {\it complex} unital $\ast$-algebra. Then
we have $\cA=\cA_h+\ii \cA_h$ be writing $a\in \cA$ as
\begin{align*}
a=a_1+\ii a_2 ,~{\rm where}~ a_1 ={\rm Re}~a:=(a^\ast {+}a)/2,~~
a_2={\rm Im}~ a :=\ii(a^\ast{-}a)/2.
\end{align*}
\subsection{Definition and Simple Properties}

Let $\cC$ be a quadratic module of  $\cA$. We denote by $\cA_b(\cC)$
the set of all elements $a \in \cA$ for which exists a a number
$\lambda_a >0$ such that
\begin{align}\label{abc}
\lambda_a{\cdot}1 \pm {\rm Re}~ a \in \cC ~{\rm and} ~\lambda_a{\cdot}1 \pm{\rm Im}~ a \in \cC.
\end{align}
The following proposition was proved in \cite{sweyl} and in \cite{cim2}.
\begin{thp}\label{p1}
{\rm (i)}  $\cA_b(\cC)$ is a unital $\ast$-subalgebra of $\cA$.\\
{\rm (ii)} An element $a\in \cA$ is in $\cA_b(\cC)$
if and only if $a^\ast a \in \cA$.
\end{thp}
We call $\cA_b(\cC)$ the {\it $\ast$-subalgebra of $\cC$-bounded
elements} of $\cA$. In the case $\cC= \sum \cA^2$ we denote
$\cA_b(\cC)$ by $\cA_b$. Note that $\cA_b(\cC)$  is the counter-part
of the ring of bounded elements (see \cite{schweig}, \cite{mm1}) in
real algebraic geometry.

 The main notion in this section is the following.

\begin{thd}
A quadratic module $\cC$ of $\cA$ is called {\rm Archimedean} if for each
element $a \in \cA_h$ there exists a $\lambda >0$ such
that $\lambda{\cdot}1-a \in \cC$ and $\lambda{\cdot}1+a \in \cC$.
\end{thd}

Let $\cC$ be a quadratic module of $\cA$. By the definition of
$\cA_b(\cC)$ the quadratic module $\cC_b{:=}\cC\cap \cA_b(\cC)$ of
the $\ast$-algebra $\cA_b(\cC)$ is Archimedean. Obviously, $\cC$ is
Archimedean if and only if $\cA_b(\cC)=\cA$. In order to prove that
a quadratic module $\cC$ is Archimedean, by Proposition \ref{p1}(i)
it suffices to show that a set of generators of $\cA$ is in
$\cA_b(\cC)$. This fact is essentially used in proving
Archimedeaness for all corresponding examples in this section.

Clearly, $\cC$ is Archimedean if and only if  $1$
is an order unit (see \cite{ja}) of the corresponding ordered vector space
$(\cA_h, \succeq)$.

Recall that a point $x$ of a subset $M$ of a real vector space $E$
is called an {\it internal point} of $M$ if for any $y \in E$ there
exists a number $\varepsilon_y >0$ such that $x+\lambda y \in M$ for
all $\lambda \in \dR$, $|\lambda| \leq \varepsilon_y$. Let $M^\circ$
denote the set of internal points of $M$.

Since order units and internal points coincide \cite{ja}, $\cC$ is Archimedean
if and only if $1$ is an internal point of $\cC$. The existence of
an internal point is the crucial assumption for Eidelheit's
separation theorem for convex sets. Let us say that a
$\ast$-representation $\pi$ is {\it $\cC$-positive} if $\pi (c) \geq
0$ for all $c \in \cC$.

\mn
\begin{thl}\label{eidel}
Let $\cC$ be an Archimedean quadratic module of $\cA$. Suppose that
$\cB$ is a convex subset of $\cA_h$ such that $\cC^\circ \cap \cB =
\emptyset$. Then there exists a state $F$ of the $\ast$-algebra
$\cA$ such that the GNS-representation $\pi_F$ is $\cC$-positive and
$F(b)\leq 0$ for all $b \in \cB$. In particular, $F$ is $\cC$-positive.
\end{thl}

\mn {\bf Proof.}  By Eidelheit's theorem (see e.g. \cite{ja}, 0.2.4)
there exists a $\dR$--linear functional $f \neq 0$ on $\cA_h$ such
that $\inf \{f(c); c \in \cC \} \geq \sup  \{f(b); b \in \cB \}$.
Because $\cC$ is a wedge, $f(c) \geq 0$ for all $c \in \cC$ and
$f(b)\leq 0$ for $b \in \cB$. Since $1 \in \cC^\circ$ and $f \neq
0$, $f(1) >0$. We extend $F:=f(1)^{-1}f$ to a $\dC$-linear
functional on $\cA$ which is denoted again by $F$. Since $
\sum~\cA^2 \subseteq \cC$, $F$ is a state of the $\ast$-algebra
$\cA$. Let $\pi_F$ denote the GNS-representation of $F$. For $c \in
\cC$ and $a \in \cA$, we have $a^\ast ca \in \cC$ and hence
$F(a^\ast ca)=f(1)^{-1} f(a^\ast ca) \geq 0$. Therefore, using
formula (\ref{gns}) we obtain
\begin{align*}
 \langle \pi_F(c)\pi_F(a)\varphi_F,\pi_F(a)\varphi_F \rangle& = \langle \pi_F(a^\ast ca) \varphi_F,\varphi_F \rangle   =
F(a^\ast ca) \geq 0,
\end{align*}
that is, $\pi_F(c) \geq 0$ and $\pi_F$ is $\cC$-positive. \hfill $\Box$
\begin{thl}\label{archbound}
If $\cC$ is an Archimedean quadratic module and $\pi$ is a
$\cC$--positive $\ast$-representation of $\cA$, then all operators
$\pi(a)$, $a \in \cA$, are bounded.
\end{thl}
\mn {\bf Proof.} Let $a \in \cA$. Since $\cC$ is Archimedean, by
Proposition \ref{p1}(ii) there exists a positive number $\lambda$ such that
$\lambda{\cdot}1 -a^\ast a \in \cC$. Therefore,
$$\langle(\pi( \lambda \cdot 1{-}a^\ast a)\varphi,\varphi \rangle = \lambda \aabs{\varphi}^2 -\aabs{\pi (a) \varphi}^2 \geq0
$$
and hence $\aabs{\pi (a) \varphi} \leq \lambda^{1/2} \aabs{\varphi}$
for all $\varphi \in \cD(\pi)$. \hfill $\Box$

\begin{thd}
A $\ast$-algebra $\cA$ is called {\rm algebraically bounded} if the
quadratic module $\sum \cA^2$ is Archimedean.
\end{thd}
Since $\ast$-representations are always $\sum \cA^2$-positive, each
$\ast$-representation of an algebraically bounded $\ast$-algebra
acts by bounded operators.

\subsection{Abstract Positivstellens\"atze for Archimedean Quadratic Modules}
For the following three propositions we assume that $\cC$ is an
Archimedean quadratic module of $\cA$.

\begin{thp}\label{ps1}
For any element $a \in \cA_h$ the following are equivalent:
\begin{itemize}
\item[\em (i)] $a+\varepsilon{\cdot}1 \in \cC$ for each
$\varepsilon >0$. \item[\em (ii)] $\pi(a) \geq 0$ for each
$\cC$--positive $\ast$-representation $\pi$ of $\cA$. \item[\em
(iii)] $f(a) \geq 0$ for each $\cC$-positive state $f$ on $\cA$.
\end{itemize}
\end{thp}
\mn {\bf Proof.} The implications (i)$\to$(ii)$\to$(iii) are
clear. To prove that (iii) implies (i) let us assume to the contrary
that $a+\varepsilon \cdot 1$ is not in $\cC$ for some $\varepsilon
>0$. Applying Lemma \ref{eidel} with $\cB:= \{a+\varepsilon\cdot 1
\}$ yields a $\cC$-positive state $f$ such that $f(a+\varepsilon
\cdot 1) \leq 0$. Then we have $f(a) <0$ which contradicts (iii). \hfill
$\Box$ \mn
\begin{thp}\label{ps2}
For $a \in \cA_h$ the following conditions are equivalent:
\begin{itemize}
\item[\em (i)] There exists $\varepsilon > 0$ such that $a-\varepsilon \cdot 1 \in
\cC$.
\item[\em (ii)] For each $\cC$-positive $\ast$-representation $\pi$ of $\cA$ there exists a number
$\delta_\pi >0$ such that $\pi(a{-}\delta_\pi {\cdot} 1) \geq0$.
\item[\em (iii)] For each $\cC$-positive state $f$ of $\cA$ there exists a number $\delta_f
>0$ such that $f(a{-}\delta_f {\cdot} 1) \geq 0$.
\end{itemize}
\end{thp}
\mn {\bf Proof.} As above, (i)$\to$(ii)$\to$(iii) is obvious.
We prove (iii)$\to$(i). Assume that (i) does not hold. We apply
Lemma \ref{eidel} to the Archimedean quadratic module $\tilde{\cC}
=\dR_+{\cdot}1 + \cC$ and $\cB=\{a\}$ and obtain a $\tilde{\cC}$-
positive state $f$ on $\cA$ such that $f(a)\leq 0$. Since $f$ is
also $\cC$-positive, this contradicts (iii). \hfill $\Box$ \mn

The assertion of next proposition I have learned  from J. Cimpric'
talk at the  Marseille conference, March 2005.
\begin{thp}\label{ps3}
For $a\in\cA_h$ the following are equivalent:
\begin{itemize}
\item[\em (i)] There exist nonzero elements $x_1, \dots, x_r$ of $\cA$ such that $\sum_{k=1}^r x_k^\ast a x_k $ belongs to $1+ \cC$.
\item[\em (ii)] For any $\cC$-positive $\ast$-representation $\pi$ of $\cA$ there
exists a vector $\eta$ such that $\langle \pi(a)\eta,\eta\rangle>0$.
\end{itemize}
\end{thp}
\mn {\bf Proof.} (i)$\to$(ii): Suppose that $\sum_k x_k^\ast a x_k
=1+c$ with $c \in  \cC$. If $\pi$ is a $\cC$-positive
$\ast$-representation and $\varphi \in \cD(\pi)$, $\varphi \neq 0$,
then
\begin{align*}
&\textstyle{\sum_k \langle \pi (a)\pi (x_k)\varphi, \pi (x_k)\varphi
\rangle} = \textstyle{\sum_k \langle \pi(x_k^\ast a x_k) \varphi,
\varphi \rangle}\\ &= \langle \pi (1+c)\varphi, \varphi \rangle \geq
\langle \pi (1) \varphi, \varphi \rangle = \aabs{\varphi}^2 >0.
\end{align*}
Hence at least one summand $\langle \pi (a)\pi (x_k)\varphi, \pi
(x_k)\varphi \rangle $ is positive.

(ii)$\to$(i): Let $\cB$ be the set of finite sums of elements
$x^\ast a x$, where $x \in \cA$, and let $\tilde{\cC}{:=} 1+\cC$. If
(i) does not hold, then $\cB \cap \tilde{\cC} = \emptyset$. By Lemma
\ref{eidel} there exists a state $f$ of $\cA$ such that the GNS
representation $\pi_f$ is $\tilde{\cC}$-positive and $f(\cB)\leq 0$.
The latter means that $f(x^\ast a x)=\langle \pi_f (a) \pi_f(x)
\varphi_f, \pi (x) \varphi_f \rangle  \leq 0$ for all $x \in \cA$.
Since $\cD(\pi_f)=\pi_f(\cA)\varphi_f$ (see e.g. \cite{unb} , 8.6),
the condition in (ii) is not satisfied for the GNS representation
$\pi_f$. \hfill $\Box$

\subsection{The Archimedean Positivstellensatz for Compact Semialgebraic Sets}

Let $f{=}(f_1,\cdots,f_k)$ be a $k$-tuple of polynomials $f_j \in
\dR[t_1,{\cdots}, t_d]$. Let $\cK_f$ be the basic closed
semialgebraic set (\ref{kf}) and $\cT_f$ the preorder  (\ref{tf})
associated with $f$. Then $\cT_f$ is a quadratic module of $\cA{=}
\dC[t_1,\dots,t_d]$. Recall that $"{\preceq}"$ denotes the order
relation defined by $\cT_f$.

\begin{thp}\label{prearchcom}
If the set $\cK_f$ is compact, then $\cT_f$ is Archimedean.
\end{thp}
\mn {\bf Proof.} Let $p \in \dR[t]$ and fix a positive number
$\lambda$ such that $\lambda^2{-}p^2> 0$ on the compact set $\cK_f$.
By Stengle's Positivstellensatz, applied to the positive polynomial
$\lambda^2{-}p^2$ on  $\cK_f$, there exist $g,h \in \cT_f$ such that
\begin{align}\label{stengle}
g(\lambda^2-p^2) =1 +h.
\end{align}
For $n \in \dN_0$, we have $p^{2n}(1+h) \in \cT_f$. Therefore, using
(\ref{stengle}) it follows that $p^{2n+2}g = p^{2n} \lambda^2 g
{-}p^{2n}(1+h) \preceq p^{2n} \lambda^2 g$. By induction we get
\begin{align}\label{p2n1}
p^{2n} g \preceq \lambda^{2n} g.
\end{align}
Since $p^{2n}(h+gp^2) \in \cT_f$, using first (\ref{stengle}) and
then (\ref{p2n1}) we obtain
\begin{align}\label{p2n2}
p^{2n} \preceq p^{2n} + p^{2n}(h+gp^2) = p^{2n}\lambda^2 g \preceq
\lambda^{2n+2}g.
\end{align}

Now we out $p:=(1+t_1^2)\cdots(1+t_d^2)$. If $|\alpha| \leq k$, $k
\in \dN$, we have
\begin{align}\label{pk}
\pm 2t^\alpha \preceq t^{2\alpha} + 1 \preceq \sum_{|\beta| \leq k} t^{2\beta} = p^k.
\end{align}
Hence there exist numbers $c >0$ and $k \in \dN$ such that $g
\preceq 2c p^k$. Combining the latter with (\ref{p2n2}), we get
$p^{2k} \preceq 2c \lambda^{2k+2} p^k$ and so
$(p^k{-}\lambda^{2k+2}c)^2 \preceq (\lambda^{2k+2}c)^2{\cdot}1$.
Therefore, by Proposition \ref{p1}(ii), $p^k{-}\lambda^{2k+2}c \in
\cA_b(\cT_f)$ and so  $p^k \in \cA_b(\cT_f)$. Since $\pm t_j \preceq
p^k$ by (\ref{pk}), we have $t_j \in \cA_b(\cT_f)$ for
$j{=}1,{\cdots},d$. From Proposition \ref{p1}(i) it follows that
$\cA_b(\cT_f)= \cA$ which means that $\cT_f$ is Archimedean. \hfill
$\Box$

\mn Using the preceding result we now give a new and {\it
elementary} proof of the author's Positivstellensatz \cite{sps}.

\begin{tht}\label{schmps}
Let $q\in \dR[t_1,{\cdots},t_d]$. If $ q(s) > 0$ for all $s \in
\cK_f$ and $\cK_f$ is compact, then $q \in \cT_f$.
\end{tht}
\mn {\bf Proof.} Assume to the contrary that $q$ is not in $\cT_f$.
By Proposition \ref{prearchcom}, $\cT_f$ is Archimedean. Therefore,
by Lemma \ref{eidel} there exists a $\cT_f$-positive state $F$ on
$\cA$ such that $F(q) \leq 0$. Let $\parallel {p} \parallel$ denote
the supremum of $p \in \dR[t]$ on the compact set $\cK_f$. Our first
aim is to show that $F$ is $\parallel \cdot \parallel $-continuous.

For let $p \in \dR[t]$. Fix $\varepsilon {>}0$ and put $\lambda{:=}
\parallel p \parallel +\varepsilon$. We define a state $F_1$ on the
polynomials in one hermitian indeterminate $x$ by
$F_1(x^n):=F(p^n)$, $n \in \dN_0$. By the solution of the Hamburger
moment problem there exists a positive Borel measure $\nu$ on $\dR$
such that $F_1(x^n)= \tint s^n d\nu(s)$, $n \in \dN_0$. For $\gamma
>\lambda$ let $\chi_\gamma$ denote the characteristic function of
$(-\infty,-\gamma]\cup[\gamma,+\infty)$. Since $\lambda^2-p^2>0$ on
$\cK_f$, we have $p^{2n} \preceq \lambda^{2n+2}g$ by equation
(\ref{p2n1}) of the preceding proof. Using the $\cT_f$-positivity of
$F$  we derive
$$
\gamma^{2n} \tint \chi_\gamma ~d\nu \leq \tint s^{2n} d\nu(s) =F_1(x^{2n})= F(p^{2n}) \leq \lambda^{2n+2} F(g)
$$
for all $n{\in} \dN$. Since $\gamma >\lambda$, the preceding implies
that $\tint \chi_\gamma~ d\nu =0$. Therefore, ${\rm supp}~\nu
\subseteq [-\lambda,\lambda]$. Using the Cauchy-Schwarz inequality
for $F$ we obtain
$$|F(p)|^2 \leq F(p^2) = F_1(x^2)=\tint_{[-\lambda,\lambda]}~ s^2 ~d\nu(s) \leq \lambda^2 = (\parallel p \parallel + \varepsilon)^2.
$$
Letting $\varepsilon \to 0$, we get $|F(p)| \leq \parallel p
\parallel$. That is, $F$ is $\parallel {\cdot} \parallel$-continuous
on $\dR[t]$.

Since $q >0$ on the compact set $\cK_f$, there is a positive number
$\delta$ such that $q-\delta \geq 0$ on $\cK_f$. By the classical
Weierstrass theorem the continuous function $\sqrt{q(s)-\delta}$  on
$\cK_f$ is uniform limit of a sequence of polynomials $p_n \in
\dR[t]$. Then $\lim_{n} \parallel p_n^2 -q +\delta \parallel =0$ and
hence $\lim_{n} F(p_n^2-q+ \delta)=0$ by the continuity of the
functional $F$. But since $F(p_n^2) \geq 0$ and $F(q)\leq 0$, we
have $F(p_n^2-q+\delta)\geq \delta >0$ which is the desired
contradiction. \hfill $\Box$

\mn
Remarks. 1. That for compact sets $\cK_f$ the preorder  $\cT_f$ is Archimedean was first shown by
T. W\"ormann \cite{wm}.\\
2. Shortly after the Positivstellensatz \cite{sps}
appeared, A. Prestel observed that there is  a small gap in the
proof. (It has to be shown that the functional $G_{n+1}$ occuring
therein is nontrivial.) This was immediately repaired by the author
and it was precisely the reasoning used in the proof of Proposition \ref{prearchcom} that filled this gap.\\
3. Having Proposition \ref{prearchcom} there are various ways to
prove Theorem \ref{schmps}. One can use the  spectral theorem as in
\cite{sps}, the Kadison-Dubois theorem as in \cite{wm}  or Jabobi's
theorem \cite{jacobi}.

\subsection{Examples of Archimedean Quadratic Modules}\label{exarch}

\mn
\begin{thex}\label{vermap} {\it Veronese Map}\\
Let $\cA$ be the complex $\ast$--algebra of rational functions generated by
$$
x_{kl}:= x_kx_l(1+x_1^2+\cdots+x_d^2)^{-1}~,~ k,l=,1,\cdots,d,
$$
where $x_0:=1$. Since $1= \sum_{r,s}  x_{rs}^2 \succeq x_{kl}^2
\succeq 0 $  for $k,l{=}1,\dots,d,$ it follows from Proposition
\ref{p1} that $x_{kl} \in \cA_b$ and hence $\cA_b=\cA$. That is, the
quadratic module $\sum \cA^2$ is Archimedean and $\cA$ is
algebraically bounded. This algebra has been used by M.Putinar and
F. Vasilescu in \cite{pv}.\hfill $\circ$
\end{thex}

A large class of algebraically bounded $\ast$--algebras is provided
by coordinate $\ast$--algebras of compact quantum groups and
quantum spaces.

\begin{thex} {\it Compact quantum group algebras}\\
Any compact quantum group algebra $\cA$ (see e.g. \cite{klimyk},
p.415) is linear span of matrix elements of finite dimensional
unitary corepresentations. These matrix elements $v_{kl}$ with
respect an orthonormal basis satisfy the relation $\sum_{l=1}^d
~v_{kl}^\ast v_{kl} =1$ for all $k$ (\cite{klimyk}, p. 401). Hence
each $v_{kl}$ is in $\cA_b$ and so $\cA_b=\cA$. That is, each
compact quantum group algebra $\cA$ is algebraically bounded  and
$\sum \cA^2$ is Archimedean.

The simplest example is  the quantum group $SU_q(2)$, $q\in \dR$.
The corresponding $\ast$--algebra  has two generators $a$ and $c$
and defining relations
$$
ac=qca,~~ c^\ast c =c c^\ast,~~ a a^\ast +q^2cc^\ast =1 ~~,~~ a^\ast a +c^\ast c=1.
$$
From the last relation we see that
$a$ and $c$ are in $\cA_b$. Hence $\cA_b=\cA$. \hfill $\circ$
\end{thex}

\begin{thex} {\it Compact quantum spaces}\\
Many compact quantum spaces have algebraically
bounded coordinate $\ast$--algebras $\cA$. Famous examples are the so-called quantum
spheres, see e.g. \cite{klimyk}, p. 449. One of the defining
relations of the $\ast$-algebra $\cA$ is
$\sum_{k=1}^n~ z_kz_k^\ast =1$ for the generators $z_1,{\dots},z_n$.
Hence we have $\cA_b=\cA$. \\
The simplest example is the  $\ast$-algebra $\cA$ with
 generators $a$ and defining relation $aa^\ast +q a^\ast
a= 1$, where $q>0$. \hfill $\circ$
\end{thex}

Weyl algebras and enveloping algebras  are not algebraically
bounded, but they do have algebraically bounded fraction
$\ast$-algebras. The fraction algebras of the next two examples have
been the main technical tools in the proofs of a strict
Positivstellensatz in \cite{sweyl} and in \cite{senv}.

\begin{thex}\label{weyl}~ {\it A fraction algebra for the Weyl algebra}\\
Let $\cW(d)$ be the Weyl algebra (Example \ref{weyl1}) and set
$N=a_1^\ast a_1+\cdots+a_d^\ast a_d$. Let us fix real number
$\alpha$ which is not an integer. Let $\cA$ be the $\ast$-subalgebra
of the fraction algebra of $\cW(d)$ generated by the elements
$$
x_{kl}:= a_ka_l(N+\alpha 1)^{-1},~k,l=0,\dots,d,~~{\rm and}~~ y_n:=(N+(\alpha+n)1)^{-1},~n \in Z,
$$
where $a_0:=1$. Then $\cA$ is algebraically bounded (\cite{sweyl}, Lemma 3.1). \hfill $\circ$
\end{thex}

\begin{thex}\label{env} {\it A fraction algebra for enveloping algebras}\\
Let $\cE(\cG)$ be the complex universal enveloping algebra of a Lie
algebra $\cG$ (Example \ref{env1}). We fix a basis
$\{x_1,\dots,x_d\}$ of the real vector space $\cG$ and put $ a:=
1+x_1^\ast x_1+ \cdots + x_d^\ast x_d$. Let $\cA$ be the unital
$\ast$--subalgebra of the fraction algebra of $\cE(\cG)$ generated
by the elements $x_{kl}:=x_k x_l a^{-1}$, $k,l{=}0,\dots,d$, where
$x_0:=1$. As shown in \cite{senv}, $\cA$ is algebraically bounded.
\hfill $\circ$
\end{thex}

\section{Transport of Quadratic Modules by Pre-Hilbert $\ast$-Bimodules}

Let $\cA$ and $\cB$ be complex unital $\ast$-algebras. We shall show
how $\cA{-}\cB$--bimodules equipped with $\cA$-- and $\cB$--valued
sesquilinear forms can be used to move quadratic modules from one
algebra to the other. Our assumptions (i)--(ix) are close to the
axioms of equivalence bimodules in the theory of Hilbert
$C^\ast$-modules (see \cite{mt}, 1.5.3).

Let $\cX$ be a left $\cA$-module and a right $\cB$-module such that
$(a{\cdot}x){\cdot}b=a{\cdot}(x{\cdot}b)$ for $a \in \cA$, $b \in
\cB$ and $x \in \cX$. Suppose that there is a  sesquilinear map
$\langle{\cdot}, {\cdot}\rangle_\cB : \cX \times \cX \to \cB$ which
is conjugate linear in the first variable and satisfies the
following conditions for $x,y\in \cX$, $b \in \cB$, and $a \in \cA$:
\begin{itemize}
\item[\rm (i)] $\langle x,y \rangle_{\cB}^\ast = \langle y,x \rangle_{\cB}$,
\item[\rm (ii)] $\langle x,y{\cdot}b \rangle_{\cB} = \langle x,y \rangle_{\cB}~ b$,
\item[\rm (iii)] $\langle a{\cdot}x, x \rangle_{\cB} = \langle x, a^\ast {\cdot}x \rangle_{\cB}$,
\item[\rm (iv)] The unit $1$ of $\cB$ is a finite sum of elements $\langle x,x \rangle_{\cB}$, where $x \in \cX$.
\end{itemize}
\mn
For a quadratic module $\cC$ of $\cA$, let $ \cC_\cX$ denote the set
of  finite sums of elements $\langle a{\cdot}x,x\rangle_{\cB} $, where
$a\in \cC$ and $x \in \cX$.
\begin{thl}\label{cbm}
$\cC_\cX$ is a quadratic module of
the $\ast$-algebra $\cB$.
\end{thl}
{\bf Proof.}  From (i) and (iii) it follows that $\cC_\cX$ is contained in
$\cB_h$. Obviously, $\cC_\cX$ is a wedge. Since the unit of $\cA$ is
in $\cC$, (iv) implies that the unit of $\cB$ is in
$\cC_\cX$. Let $b \in \cB$, $c\in \cC$ and $x \in \cX$. Using
conditions (ii) and (i) and the bimodule axiom
$(a{\cdot}x){\cdot}b= a{\cdot}(x{\cdot}b)$ we obtain
\begin{align*}
&b^\ast \langle a{\cdot}x,x \rangle_{\cB} b= b^\ast \langle a{\cdot}x,x {\cdot}b \rangle_{\cB}
= (\langle x{\cdot}b, a{\cdot}x \rangle_{\cB} b)^\ast =\\
&(\langle x{\cdot}b,(a{\cdot}x){\cdot}b \rangle_{\cB} )^\ast =\langle x{\cdot}b,a{\cdot}(x{\cdot}b) \rangle_{\cB} )^\ast
= \langle a{\cdot}(x{\cdot}b),(x{\cdot}b)\rangle_{\cB}.
\end{align*}
Hence $\cC_\cX$ satisfies (\ref{QM2}) and $\cC_\cX$ is a quadratic
module. \hfill $\Box$

\mn Suppose that $\langle{\cdot}, {\cdot}\rangle_\cA : \cX \times
\cX \to \cA$ is a sesquilinear map which is conjugate linear in the
second variable such that for $x,y\in \cX$, $b \in \cB$ and $a \in
\cA$:
\begin{itemize}
\item[\rm (v)] $\langle x,y \rangle_{\cA}^\ast = \langle y,x \rangle_{\cA}$,
\item[\rm (vi)] $\langle a{\cdot}x,y \rangle_{\cA} = a~\langle x,y \rangle_{\cA}$,
\item[\rm (vii)] $\langle x{\cdot}b, x \rangle_{\cA} = \langle x, x{\cdot}b^\ast \rangle_{\cA}$,
\item[\rm (viii)] The unit $1$ of $\cA$ is a finite sum of elements $\langle x,x \rangle_{\cA}$, where $x \in \cX$.
\end{itemize}
\mn For a quadratic module $\cP$ of $\cB$, let $_{\cX}\cP $ denote
the finite sums of elements $\langle y,y{\cdot}b\rangle_{\cA} $,
where $b\in \cP$ and $y \in \cX$. A similar reasoning as in the
proof of Lemma \ref{cbm} shows {\it $_{\cX}\cP$ is a quadratic
module for $\cA$}.

Finally, we  assume the following compatibility condition:
\begin{itemize}
\item[\rm (ix)]~~ $\langle x,y \rangle_\cA \cdot z = x\cdot \langle y,z \rangle_\cB$ for all $x,y,z \in \cX.$
\end{itemize}
\begin{thp}
%Suppose that all axioms (i)--(ix) are fulfilled.
If $\cC$ is a quadratic module of $\cA$ and $\cP$ is a quadratic
module of $\cB$, then we have ${_\cX}(\cC_\cX) \subseteq \cC$ and
$(_{\cX}\cP)_\cX \subseteq \cP$.
\end{thp}
\mn {\bf Proof.} We prove only the first inclusion. Since
${_\cX}(\cC_\cX)$ consists of sums of elements of the form $\langle
y,y{\cdot} {\langle a{\cdot}x,x \rangle_\cB} \rangle_\cA$, where $a
\in \cC$ and $x,y \in \cX$, it suffices to show that these elements
are in $\cC$. We compute
\begin{align*}
\langle y,y{\cdot} {\langle a{\cdot}x,x \rangle_\cB} \rangle_\cA &=
\langle y {\cdot}{\langle x,a{\cdot}x \rangle_\cB}, y \rangle_\cA =
\langle \langle y,x\rangle_\cA {\cdot}(a{\cdot}x),y \rangle_\cA \\ &=
\langle y,x \rangle_\cA \langle a{\cdot}x,y \rangle_\cA=
({\langle x,y\rangle_\cA})^\ast ~a~ \langle x,y\rangle_\cA,
\end{align*}
where the first equality follows from assumptions (vii) and (i), the
second from (ix), the third from (vi), and the fourth from (v) and
(vi). By (\ref{QM2}) the terms on the right hand side of the
preceding equations is in $\cC$. \hfill $\Box$

\mn
We illustrate these general constructions by an important example.

\begin{thex}
{\it Quadratic modules of~ k-positive $n \times n$ matrices}\\
Let $R$ be a unital $\ast$-algebra. Set $\cA{=}\cM_{k}(R)$,
$\cB{=}\cM_{n}(R)$, and $\cX{=}\cM_{kn}(R)$. Then $\cX$ is an
$\cA{-}\cB$--bimodule with module operations defined by the left
resp. right multiplications of matrices and equipped with
$\cB$--resp. $\cA$-valued "scalar products" $\langle x,y
\rangle_{\cB}:=x^\ast y$ and $\langle x,y \rangle_{\cA}:=x y^\ast$
for $x,y \in \cX$. With these definitions all assumptions (i)--(ix)
are satisfied.

If $\cC$ and $\cP$ are quadratic modules for $\cA{=}\cM_{k}(R)$ and
$\cB{=}\cM_n(R)$, respectively,, then the  quadratic module
$\cC_\cX$ and $_{\cX}\cP$ are given by
\begin{align}\label{cnk1}
&\cC_{n,k}:=\cC_\cX= \left\{ \textstyle{\sum_{l=1}^s ~ x_l^\ast a_l
x_l};~~ a_l \in \cC,~ x_l \in  \cM_{kn}(R),~ s \in \dN ~ \right\},\\
\nonumber & \cP_{k,n}:={_{\cX}\cP}= \left\{ \textstyle{\sum_{l=1}^s
~ y_l b_l y_l^\ast};~~ b_l \in \cP,~ y_l \in  \cM_{kn}(R),~ s \in
\dN ~ \right\}.
\end{align}

Now we specialize the preceding by setting
$R{=}\dC[t_1,{\dots},t_d]$. Let  $\cC$ be  the set $\cM_k(\dC[t])_+$
of hermitian $k \times k$ matrices over $\dC[t]$ which are positive
semidefinite for all $s \in \dR^d$. Put $\cC_{n,0}:= \sum \cB^2$.
From (\ref{cnk1})  we obtain an increasing chain of quadratic
modules
\begin{align}\label{cnk2}
\cC_{n,0} \subseteq \cC_{n,1} \subseteq \cC_{n,2}\subseteq \dots \subseteq \cC_{n,n}
\end{align}
of $\cB{=}\cM_{n}(\dC[t])$. Matrices belonging to $\cC_{n,k}$ will be called {\it $k$-positive}.

If $d=1$ and $a \in \cC_{n,n}$, then the matrix $a$ is positive
semidefinite on $\dR$ and hence of the form $a=b^\ast b$ for some $b
\in \cM_n(\dC[t])$ \cite{dj}. Therefore all quadratic modules in (\ref{cnk2}) coincide
with $\sum \cM_{n}(\dC[t])^2$.

Suppose now that $d\geq 2$. As
shown in \cite{frs}, the matrix
$$
\left(
  \begin{array}{ccc}
    1+t_1^4t_2^2 & t_1t_2 \\
    t_1t_2 & 1+t_1^2t_2^4 \\
  \end{array}
\right).
$$
is in $\cC_{2,2}$, but not in $\cC_{2,1}$. For $d\geq 2$ we have a sequence $\cC_{n,k}$ of
intermediate quadratic modules between the two extremes
$\cC_{n,0}=\sum \cM_{n}(\dC[t])^2$ and $\cC_{n,n}= \cM_{n}(\dC[t])_+$. These
quadratic modules are used in Hilbert space representation theory to
characterize $k$-positive representations of the polynomial algebra
$\dC[t_1,{\dots},t_d]$  (see \cite{frs} and \cite{unb}, Proposition
11.2.5).
\end{thex}


\begin{thebibliography}{10}
\bibitem{bcr}
{\sc J. Bochnak, M. Coste and M.-F. Roy}, {\sl Real Algebraic
Geometry}, Springer-Verlag, Berlin, 1998.

\bibitem{cim1}
{\sc J. Cimpric}, {\sl Maximal quadratic modules on $\ast$-rings},
2005, to appear in Algebras and Representation theory.

\bibitem{cim2}
{\sc J. Cimpric}, {\sl A representation theorem for quadratic
modules on $\ast$-rings}, 2005, to appear in Canadian Mathematical
Bulletin.

\bibitem{dj}
{\sc D. Z. Djokovic}, {\sl Hermitean matrices over polynomial
rings}, J. Algebra {\bf 43} (1976), pp. 359--374.

\bibitem{frs}
{\sc J. Friedrich and K. Schm\"udgen}, {\sl n-Positivity of
unbounded $\ast$-representations}, Math. Nachr. {\bf 141} (1989),
pp. 233--250.

\bibitem{gam}
{\sc F.R. Gantmacher}, {\sl Matrizentheorie}, DVW, Berlin, 1986.

\bibitem{gr}
{\sc D. Gondard and P. Ribenboim}, {\sl Le 17e probleme de Hilbert
pour les matrices}, Bull. Sci. Math. {\bf 98} (1974), pp. 49--56.

\bibitem{hel}
{\sc J.W. Helton}, {\sl Positive noncommutative polynomials are sums
of squares}, Ann. Math. {\bf 156} (2002), pp. 675--694.

\bibitem{hm}
{\sc J.W. Helton and S. McCullough}, {\sl A Positivstellensatz for
non-commutative polynomials}, Trans. Amer. Math. soc. {\bf 356}
(2004), pp. 3721--3737.

\bibitem{hmp}
{\sc J.W. Helton, S. McCullough and M. Putinar}, {\sl A
non-commutative Positivstellensatz on isometries}, J. reine angew.
Math. {\bf 568} (2004), pp. 71--80.

\bibitem{jacobi}
{\sc T. Jacobi}, {\sl A representation theorem for certain partially
ordered commutative rings}, Math. Z. {\bf 237} (2001), pp. 259--273.

\bibitem{ja}
{\sc G. Jameson}, {\sl Ordered Linear Spaces}, Lecture Notes in
Math. No. {\bf 141}, Spinger-Verlag, Berlin, 1970.

\bibitem{klimyk}
{\sc A. Klimyk and K. Schm\"udgen}, {\sl Quantum Groups and Their
Representations}, Springer-Verlag, Berlin, 1997.

\bibitem{mt}
{\sc V.M. Manuilov and E.V. Troitsky}, {\sl Hilbert
$C^\ast$-Modules}, Amer. Math. Soc., Transl. Math. Monographs {\bf
226}, 2005.

\bibitem{mm}
{\sc M. Marshall}, {\sl Positive Polynomials and sums of Squares},
Univ. Pisa, Dipart. Mat. Istituti Editoriali e Poligrafici
Internaz., 2000.

\bibitem{mm1}
{\sc M. Marshall}, {\sl Extending the Archimedean Positivstellensatz
to the non-compact case}, Can. Math. Bull. {\bf 14}(2001), 223--230.

\bibitem{os}
{\sc V. Ostrovskyj and Yu. Samoilenko}, {\sl Introduction to the
Theory of Finitely Presented $\ast$-Algebras}, Harwood Acad. Publ.,
1999.

\bibitem{pd}
{\sc A. Prestel and Ch. N. Delzell}, {\sl Positive Polynomials},
Spinger-Verlag, Berlin, 2001.

\bibitem{ps}
{\sc C. Procesi and M. Schacher}, {\sl A Non-Commutative Real
Nullstellensatz and Hilbert's 17th Problem}, Ann Math. {\bf 104}
(1976), pp. 395--406.

\bibitem{pv}
{\sc M. Putinar and F. Vasilescu}, {\sl Solving moment problems by
dimension extension}, Ann. Math. {\bf 149}, pp. 1087--1107.

\bibitem{loc}
{\sc H. Sch\"afer}, {\sl Topological Vector Spaces},
Springer-Verlag, Berlin,1972

\bibitem{rev}
{\sc K. Schm\"udgen}, {\sl Graded and filtrated topological
$\ast$-algebras II. The closure of the positive cone}, Rev. Roum.
Math. Pures et Appl. {\bf 29} (1984), pp. 89--96.

\bibitem{unb}
{\sc K. Schm\"udgen}, {\sl Unbounded Operator Algebras and
Representation Theory}, Birkh\"auser-Verlag, Basel, 1990.

\bibitem{sps}
{\sc K. Schm\"udgen}, {\sl The $K$-moment problem for compact
semi-algebraic sets}, Math. Ann. {\bf 289} (1991), pp. 203--206.

\bibitem{sweyl}
{\sc K. Schm\"udgen}, {\sl A strict Positivstellensatz for the Weyl
algebra}, Math. Ann. {\bf 331} (2005), pp. 779--794.

\bibitem{senv}
{\sc K. Schm\"udgen}, {\sl A strict Positivstellensatz for
enveloping algebras}, Math. Z. {\bf 254} (2006), pp. 641--653.

\bibitem{schweig}
{\sc M. Schweighofer}, {\sl Iterated rings of bounded elements and
generalizations of Schm\"udgen's theorem}, Dissertation, Konstanz,
2002.

\bibitem{wm}
{\sc T. W\"ormann}, {\sl Strict positive Polynome in der
semialgebraischen Geometrie}, Dissertion, Universit\"at Dortmund,
1998.

\end{thebibliography}
\end{document}